\documentclass[12pt]{article}%
\usepackage{amsmath,amssymb,amsfonts,amsthm}
\usepackage{latexsym}
\usepackage[english]{babel}
\usepackage{showidx}
\usepackage{graphicx}
\usepackage{graphics}
\usepackage{babel,amsmath,amssymb,amsbsy,amsfonts,latexsym,amsthm,mathrsfs}

\def\I{\mbox{\large \bf 1}}
\def\P{{\mathbb P}}
\def\R{{\mathbb R}}
\def\E{{\mathbb E}}
\def\D{{\mathbb D}}
\def\N{{\mathbb N}}
\def\cl#1{{\mathscr #1}}

\newtheorem{theorem}{Theorem}

\newtheorem{definition}[theorem]{Definition}

\newtheorem{lemma}[theorem]{Lemma}

\newtheorem{proposition}[theorem]{Proposition}
\newtheorem{remark}[theorem]{Remark}

\begin{document}

\parindent 0pt

\title{
\textbf{Riesz transform and
integration by parts formulas for random variables
}}
\author{
{\sc Vlad Bally}\footnote{Laboratoire d'Analyse et de Math\'ematiques Appliqu\'ees, UMR 8050,
Universit\'e Paris-Est Marne-la-Vall\'ee, 5 Bld
Descartes, Champs-sur-Marne, 77454 Marne-la-Vall\'ee Cedex 2, France.
Email:
\texttt{bally@univ-mlv.fr}
}\smallskip\\
{\sc Lucia Caramellino}\footnote{Dipartimento di Matematica, Universit\`a di Roma - Tor Vergata,
Via della Ricerca Scientifica 1, I-00133 Roma, Italy. Email: \texttt{caramell@mat.uniroma2.it}}\smallskip\\
}
\date{}
\maketitle

{\textbf{Abstract.}} 
We use integration by parts formulas to give estimates for the $L^p$ norm
of the Riesz transform. This is motivated by the representation formula
for conditional expectations of functionals on the Wiener space already given
in Malliavin and Thalmaier \cite{bib:[M.T]}.
As a consequence, we obtain  regularity and estimates for 
the density of non degenerated functionals on the Wiener space. 
We also give a semi-distance which characterizes the convergence to
the boundary of the set of the strict positivity points for the density.

\medskip


\medskip


\section{Introduction}

The starting point of this paper is the representation theorem for densities
and conditional expectations of random variables based on the
Riezs transform, recently given by Malliavin and Thalmaier in \cite{bib:[M.T]}. Let us recall it.
Let $F$ and $G$ denote  random variables taking values on $\R^d$ and $\R$ respectively
and consider the following integration by parts
formula: there exist some integrable random variables $H_{i}(F,G)$ such
that for every test function $f\in C_{c}^{\infty }(\R^{d})$
\[
IP_{i}(F,G)\quad \E(\partial _{i}f(F)G)=-\E(f(F)H_{i}(F,G)),\quad i=1,...,d.
\]%
Malliavin and Thalmaier proved that if 
$%
IP_{i}(F,1),i=1,...d$ hold and the law of $F$ has a continuous density $%
p_{F}$, then
\[
p_{F}(x)=-\sum_{i=1}^{d}\E(\partial _{i}Q_{d}(F-x)H_{i}(F,1))
\]%
where $Q_{d}$ denotes the Poisson kernel on $\R^{d}$ (that is the fundamental
solution of the equation $\delta _{0}=\Delta Q_{d}).$ Moreover, they proved
also that if $IP_{i}(F,G),i=1,...d$,  a similar
representation formula holds also for the conditional expectation of $G$
with respect to $F$. The interest of Malliavin and Thalmaier in this
representations come from numerical reasons - this allows one to simplify the
computation of densities and conditional expectations using a
Monte Carlo method. This is crucial in order to implement numerical
algorithms for solving non linear PDE's or optimal stopping
problems - for example for pricing American options. But there is a
difficulty coming on here: the variance of the estimators produced by such a
representation formula is infinite. Roughly speaking, this comes from the
blowing up of the Poisson kernel around zero: $\partial _{i}Q_{d}\in L^{p}$ for
$p=d/(d-1)<2$ but not for $p=2$. So estimates of $\E(\left\vert \partial
_{i}Q_{d}(F-x)\right\vert ^{p})$ are crucial in this framework and this is
the central point of interest in our paper. In
\cite{bib:[K-H1]} and \cite{bib:[K-H2]},
Kohasu-Higa and Yasuda 
proposed a solution to this problem using some cut off arguments. And in
order to find the optimal cut off level they used the estimates of $%
\E(\left\vert \partial _{i}Q_{d}(F-x)\right\vert ^{p})$ which we prove in
this paper (actually, they used a former version given in the preprint \cite{bib:[B.C]}).

\smallskip

So our central result concerns estimates of $\E(\left\vert \partial
_{i}Q_{d}(F-x)\right\vert ^{p}).$ It turns out that, in addition to the interest in
numerical problems, such estimates represent a suitable instrument in order
to obtain regularity of the density of functionals on the Wiener space - for
which Malliavin calculus produces integration by parts formulas. 
 Before going further let us mention that one may also
consider integration by parts formulas of higher order, that is
\[
IP_{\alpha }(F,G)\quad E(\partial _{\alpha }f(F))=E(f(F)H_{\alpha }(F,G))
\]%
where $\alpha =(\alpha _{1},...,\alpha _{k}).$ We say that 
an integration by parts formula of order $k$ holds if
this is true for every $\alpha\in \{1,\ldots,d\}^k$. 
 Now, a first question is: which
is the order $k$ of integration by parts that one needs in order to prove
that the law of $F$ has a continuous density $p_{F}?$ If one employs a
Fourier transform argument (see Nualart \cite{bib:[N]}) or the representation of the
density by means of the Dirac function (see Bally \cite{bib:[B]}) then one needs $d$
integration by parts if $F\in \R^{d}.$ In \cite{bib:[M1]} Malliavin proves that
integration by parts of order one is sufficient in order to obtain a
continuous density, does not matter the dimension $d$ (he employs some
harmonic analysis arguments). A second problem concerns estimates of the
density $p_{F}$ (and of its derivatives) and such estimates involve the $%
L^{p}$ norms of the weights $H_{\alpha }(F,1).$ In the approach using the Fourier
transform or the Dirac function, $\Vert H_{\alpha }(F,1)\Vert
_{p},\vert \alpha \vert \leq d$ are involved if one estimates $%
\Vert p_{F}\Vert _{\infty }$. But in \cite{bib:[Sh]} Shigekawa obtains
estimates of $\Vert p_{F}\Vert _{\infty }$ depending only on $%
\Vert H_{i}(F,1)\Vert _{p},$ so on the weights of order one (and similarly for
derivatives). In order to do it, he needs some Sobolev inequalities
that he proves using a representation formula based on the Green
function and some estimates of modified Bessel functions. Our program and
our results are similar but the instrument used in our paper is the Riesz
transform and the estimates of the Poisson kernel mentioned above.

\smallskip

Let us be more precise. Notice that $IP_{i}(F,G)$\ may also be written as
\[
IP_{i}(F,G)\quad \int \nabla f(x)g(x)\mu _{F}(dx)=\int f(x)\partial ^{\mu
_{F}}g(x)\mu _{F}(dx)
\]%
where $\mu _{F}$ is the law of $F,g(x):={\mathbb{E}}(G\mid F=x)$ and $%
\partial ^{\mu _{F}}g(x):={\mathbb{E}}(H(F,G)\mid F=x).$ This suggests that
we can work in the abstract framework of Sobolev spaces with respect to the
probability measure $\mu _{F}$ (instead of the usual Lebesgue measure). More
precisely for a probability measure $\mu $ we denote by $W_{\mu }^{1,p}$ the
space of the functions $\phi \in L^{p}(d\mu )$ for which there exists some
functions $\theta _{i}\in L^{p}(d\mu ),i=1,\ldots ,d$ such that $\int \phi
\partial _{i}fd\mu =-\int \theta _{i}fd\mu $ for every test function $f\in
C_{c}^{\infty }({\mathbb{R}}^{d}).$ If $F$ is a random variable of law $\mu $
then the above duality relation reads ${\mathbb{E}}(\phi (F)\partial
_{i}f(F))=-{\mathbb{E}}(\theta _{i}(F)f(F))$ and these are the usual
integration by parts formulas in the probabilistic approach - for example $%
-\theta _{i}(F)$ is connected to the weight produced by Malliavin calculus for a
functional $F$ on the Wiener space. But one may consider other frameworks -
as the Malliavin calculus on the Wiener Poisson space for example. This
formalism has already been used by Shigekawa in \cite{bib:[Sh]} and a slight variant
appears in the book of Malliavin and Thalmaier \cite{bib:[M.T]} (the so called covering
vector fields).

In Section \ref{sect-sobolev} we prove the following result: if $1\in W_{\mu
}^{1,p}$ (or equivalently $IP_i(F,1), i=1,\ldots,d$ hold) for some $p>d$ then
\[
\sup_{x\in \R^{d}}\sum_{i=1}^{d}\int \left\vert \partial
_{i}Q_{d}(y-x)\right\vert ^{p/(p-1)}\mu (dy)\leq C_{d,p}\left\Vert
1\right\Vert _{W_{\mu }^{1,p}}^{\ell _{d,p}},
\]%
Moreover $\mu (dx)=p_{\mu }(dx)dx
$, with $p_{\mu }$ H\"{o}lder continuous of order $1-d/p$, and the following
representation formula holds:
\[
p_{\mu }(x)=\sum_{i=1}^{d}\int \partial _{i}Q_{d}(y-x)\partial _{i}^{\mu
}1(y)d\mu (y).
\]%
%
More generally, let $\mu _{\phi }(dx):=\phi (x)\mu (dx).$ If $\phi \in
W_{\mu }^{1,p}$ then $\mu _{\phi }(dx)=p_{\mu _{\phi }}(x)dx$ and $p_{\mu
_{\phi }}$ is H\"{o}lder continuous. This last generalization is important
from a probabilistic point of view because it produces a representation
formula and regularity properties for the conditional expectation. We
introduce in a straightforward way higher order Sobolev spaces $W_{\mu
}^{m,p},m\in \N$ and we prove that if $1\in W_{\mu }^{m,p}$ then $p_{\mu }$
is $m-1$ times differentiable and the derivatives of order $m-1$ are H\"{o}%
lder continuous. And the analogous result holds for $\phi \in W_{\mu }^{m,p}.
$ So if we are able to iterate $m$ times the integration by parts formula we
obtain a density in $C^{m-1}.$ 
These results are already obtained by
Shigekawa. In our paper we get some supplementary information about the
queues of the density function and we develop more the applications to
conditional expectations.

Furthermore, we prove an alternative representation formula. Suppose that $F$
satisfies integration by parts formulas in order to get that its law $\mu$
has a $C^1$ density $p_\mu$. We set $U_\mu=\{p_\mu>0\}$ and for $x,y\in U_\mu
$, $A_{x,y}=\{\varphi\,:\,[0,1]\to U_\mu\,;\, \varphi\in C^1,\varphi_0=x,
\varphi_1=y\}$. Then for any $\varphi\in A_{x,y}$ one has
\[
p_\mu(y)=p_\mu(x)\exp\Big(\int_0^1\langle \partial
^\mu1(\varphi_t),\dot\varphi_t\rangle dt\Big),
\]
a formula which generalizes the one given by Bell \cite{bib:bell} (he
assumes $U_\mu={\mathbb{R}}^d$ and takes $\varphi$ as the straight line).
The above formula suggests to introduce the following Riemannian
semi-distance on $U_\mu$: setting $C^{ij}_\mu=\partial^\mu_i 1
\partial^\mu_j 1$, $i,j=1,\ldots,d$, for $x,y\in U_\mu$ one defines
\[
d_\mu(x,y)=\inf\Big\{\int_0^1\langle C_\mu(\varphi_t)\dot\varphi_t,
\dot\varphi_t\rangle^{1/2} dt\,;\, \varphi\in A_{x,y}\Big\}.
\]
Such a distance is of interest in the following framework. Suppose that $F$
is a non degenerated and smooth r.v. on the Wiener space, so that
integration by parts formulas hold and $F$ has a smooth probability density.
If $F$ is one dimensional then Fang \cite{bib:[F]} proved that $\overline{U}%
_\mu$ is connected and the interior of $\overline{U}_\mu$ is given by $U_\mu$
(so the density is itself strictly positive in the interior of the support
of the law). But this is false in the multidimensional case, as shown by D.
Nualart \cite{bib:[N]} through a counterexample. \label{mall-conj} Then
Malliavin suggested that one has to replace the Euclidean distance by the
intrinsic distance associated to the Dirichlet form of $F$ (see Hirsch and
Song \cite{bib:[H.S]} for details). And he conjectured that if $d(\cdot,\cdot)$ is such a
distance then for any sequence $\{x_n\}_n\subset U_\mu$ such that $%
p_\mu(x_n)\to 0$ then $d(x_n,x_1)\to\infty$, as $n\to\infty$. But Hirsch and
Song \cite{bib:[H.S]} provided a counterexample which shows that it is false
as long as the intrinsic distance is taken into account. We prove here that
the Malliavin's conjecture is true but with the distance $d_\mu$ defined
above (and in fact, we prove the equivalence).

\smallskip

The paper is organized as follows. In Section \ref{sect-sobolev} we develop
the main results in the abstract Sobolev spaces framework. In Section \ref%
{sect-IBP} we translate these results in probabilistic terms and in Section %
\ref{sect-acc} we give their applications on the Wiener space.

\section{Sobolev spaces associated to a probability measure and Riesz
transform}\label{sect-sobolev}

\subsection{Definitions and main objects}
We consider a probability measure $\mu $ on $\R^{d}$ (with the Borel $\sigma $%
-field) and we denote by $L_{\mu }^{p}=\{\phi :\int \left\vert \phi
(x)\right\vert ^{p}\mu (dx)<\infty \}$ and we put $\left\Vert \phi
\right\Vert _{L_{\mu }^{p}}=(\int \left\vert \phi (x)\right\vert ^{p}\mu (dx))^{1/p}.$ We also denote by $W_{\mu }^{1,p}$ the space of the functions $%
\phi \in L_{\mu }^{p}$ for which there exists some functions $\theta _{i}\in
L_{\mu }^{p},i=1,\ldots ,d$ such that, for every test function $f\in
C_{c}^{\infty }(\R^{d}),$ one has
\begin{equation*}
\int \partial _{i}f(x)\phi (x)\mu (dx)=-\int f(x)\theta _{i}(x)\mu (dx).
\end{equation*}%
We denote $\partial _{i}^{\mu }\phi =\theta _{i}.$ And we define the norm%
\begin{equation*}
\left\Vert \phi \right\Vert _{W_{\mu }^{1,p}}=\left\Vert \phi \right\Vert
_{L_{\mu }^{p}}+\sum_{i=1}^{d}\left\Vert \partial _{i}^{\mu }\phi
\right\Vert _{L_{\mu }^{p}}.
\end{equation*}

We similarly define the Sobolev spaces of higher order. Let $%
\alpha =(\alpha _{1},\ldots ,\alpha _{m})\in \{1,\ldots ,d\}^{m}$ be a multi index.
We denote by $\left\vert \alpha \right\vert =m$ the length of $\alpha $ and
for a function $f\in C^{m}(\R^{d})$ we denote by $\partial _{\alpha
}f=\partial _{x_{\alpha _{1}}}\ldots \partial _{x_{\alpha _{1}}}f$ the standard
derivative corresponding to the multi index $\alpha .$ Then we define $%
W_{\mu }^{m,p}$ to be the space of the functions $\phi \in L_{\mu }^{p}$
such that for every multi index $\alpha $ with $\left\vert \alpha
\right\vert \leq m$ there exists some functions $\theta _{\alpha }\in L_{\mu
}^{p}$ such that%
\begin{equation*}
\int \partial _{\alpha }f(x)\phi (x)\mu (dx)=(-1)^{\left\vert \alpha
\right\vert }\int f(x)\theta _{\alpha }(x)\mu (dx)\quad \forall f\in
C_{c}^{\infty }(\R^{d}).
\end{equation*}%
We denote $\partial _{\alpha }^{\mu }\phi =\theta _{\alpha }$ and we define
the norm%
\begin{equation*}
\left\Vert \phi \right\Vert _{W_{\mu }^{m,p}}=\left\Vert \phi \right\Vert
_{L_{\mu }^{p}}+\sum_{\left\vert \alpha \right\vert \leq m}\left\Vert
\partial _{\alpha }^{\mu }\phi \right\Vert _{L_{\mu }^{p}}.
\end{equation*}%
We will use the notation $L^{p},W^{m,p}$ for the spaces associated to the
Lebesgue measure (instead of $\mu )$, which are the standard $L^{p}$ and the
standard Sobolev spaces which are used in the literature. If $D\subset \R^{d}$
is an open set we denote by $W_{\mu }^{m,p}(D)$ the space of the functions $%
\phi $ which verify the integration by parts formula for test functions $f$
which have a compact support included in $D$ (so $W_{\mu }^{m,p}=W_{\mu
}^{m,p}(\R^{d})).$ The same for $W^{m,p}(D),L^{p}(D),L_{\mu }^{p}(D).$

Our aim is to study the link between $W_{\mu }^{m,p}$ and $W^{m,p}$ and the
main tool is the Riesz transform that we introduce now. The fundamental
solution $Q_{d}$ of the equation $\Delta Q_{d}=\delta _{0}$\ in $\R^{d}$ has
the following explicit form:
\begin{equation}
Q_{2}(x)=a_{2}^{-1}\ln \left\vert x\right\vert \quad \mbox{and}\quad
Q_{d}(x)=-a_{d}^{-1}\left\vert x\right\vert ^{2-d},d>2  \label{den4}
\end{equation}%
where $a_{d}$ is the area of the unit sphere in $\R^{d}.$ For $f\in
C_{c}^{1}(\R^{d})$ one has
\begin{equation*}
f=(\nabla Q_{d})\ast \nabla f.
\end{equation*}%
In Theorem 4.22 of Malliavin and Thalmaier \cite{bib:[M.T]}, this representation for the function $f$ is called  the
Riesz transform of $f$ and is employed in order to obtain representation
formulas for the conditional expectation. Moreover, some analogues
representation formulas for functions on the sphere and on the ball are used
by Malliavin and E. Nualart in \cite{bib:[M.N]} in order to obtain lower bounds for
the density of a strongly non degenerated random variable.\

Setting $A_{2}=1$ and $A_{d}=d-2$, we have
\begin{equation}
\partial _{i}Q_{d}(x)=a_{d}^{-1}A_{d}\frac{x_{i}}{\left\vert x\right\vert
^{d}}.  \label{den5}
\end{equation}%
By using polar coordinates, one has
\begin{equation}
\int_{\left\vert x\right\vert \leq 1}\left\vert \partial
_{i}Q_{d}(x)\right\vert ^{1+\delta }dx\leq A_{d}^{1+\delta
}\int_{0}^{1}\left\vert \frac{r}{r^{d}}\right\vert ^{1+\delta
}r^{d-1}dr=A_{d}^{1+\delta }\int_{0}^{1}\frac{1}{r^{\delta (d-1)}}dr
\label{den6}
\end{equation}%
which is finite for any $\delta <\frac{1}{d-1}.$ But $\left\vert \partial
_{i}^{2}Q_{d}(x)\right\vert \sim \left\vert x\right\vert ^{-d}$ and so $%
\int_{\left\vert x\right\vert \leq 1}\left\vert \partial
_{i}^{2}Q_{d}(x)\right\vert dx=\infty .$ This is the reason for which we
have to integrate by parts once and to remove one derivative, but we may
keep the other derivative.

In order to include the one dimensional case we set $Q_{1}(x)=\max
\{x,0\},a_{1}=A_{1}=1$ and we have
\begin{equation*}
\frac{dQ_{1}(x)}{dx}=1_{(0,\infty )}(x).
\end{equation*}%
In this case the above integral is finite for every $\delta >0.$

\subsection{An absolute continuity criterion}
For a function $\phi \in L_{\mu }^{1}$ we denote by $\mu _{\phi }$ the
signed finite measure defined by
\begin{equation*}
\mu _{\phi }(dx):=\phi (x)\mu (dx).
\end{equation*}%

We prove now the following theorem, which is
starting point of our next results.

\begin{theorem}\label{th1}
\textbf{A.}
Let $\phi \in W_{\mu }^{1,1}.$ Then%
\begin{equation*}
\int \left\vert \partial _{i}Q_{d}(y-x)\partial _{i}^{\mu }\phi
(y)\right\vert \mu (dy)<\infty
\end{equation*}%
for a.e. $x\in \R^d$ and $\mu _{\phi }(dx)=p_{\mu _{\phi }}(x)dx$ with%
\begin{equation}
p_{\mu _{\phi }}(x)=-\sum_{i=1}^{d}\int \partial _{i}Q_{d}(y-x)\partial
_{i}^{\mu }\phi (y)\mu (dy).  \label{rep0}
\end{equation}%

\textbf{B.}
If $\phi \in W_{\mu }^{m,p},p\geq 1$ for some $m\geq 2$,
then %
\begin{equation}
\partial _{\alpha }p_{\mu _{\phi }}(x)=-\sum_{i=1}^{d}\int \partial
_{i}Q_{d}(y-x)\partial _{(\alpha ,i)}^{\mu }\phi (y)\mu (dy)  \label{rep2}
\end{equation}%
where $\alpha $ is any multi index of length less or equal to $m-1$.
If in addition $1\in W_{\mu }^{1,p}$, the following alternative
representation formula holds:
\begin{equation}
\partial _{\alpha }p_{\mu _{\phi }}(x)=p_{\mu }(x)\partial _{\alpha }^{\mu
}\phi (x).  \label{rep3}
\end{equation}%
In particular, taking $\phi =1$ and $\alpha =\{i\}$ one has
\begin{equation}
\partial _{i}^{\mu }1=1_{\{p_{\mu }>0\}}\partial _{i}\ln p_{\mu }.
\label{rep6}
\end{equation}
\end{theorem}

\textbf{Proof}. \textbf{A.}
We take $f\in C_{c}^{1}(\R^{d})$ and we write $f=\Delta
(Q_{d}\ast f)=\sum_{i=1}^{d}(\partial _{i}Q_{d})\ast (\partial _{i}f).$ Then
\begin{eqnarray*}
\int fd\mu _{\phi } &=&\int f\phi d\mu
=\sum_{i=1}^{d}\int \mu (dx)\phi (x)\int \partial _{i}Q_{d}(z)\partial
_{i}f(x-z)dz\\
&=&\sum_{i=1}^{d}\int \partial _{i}Q_{d}(z)\int \mu (dx)\phi
(x)\partial _{i}f(x-z))dz \\
&=&-\sum_{i=1}^{d}\int \partial _{i}Q_{d}(z)\int \mu (dx)f(x-z)\partial
_{i}^{\mu }\phi (x)dz\\
&=&-\sum_{i=1}^{d}\int \mu (dx)\partial _{i}^{\mu }\phi
(x)\int \partial _{i}Q_{d}(z)f(x-z)dz \\
&=&-\sum_{i=1}^{d}\int \mu (dx)\partial _{i}^{\mu }\phi (x)\int \partial
_{i}Q_{d}(x-y)f(y)dy \\
&=&\int f(y)\Big(-\sum_{i=1}^{d}\int \partial _{i}Q_{d}(x-y)\partial _{i}^{\mu
}\phi (x)\mu (dx)\Big)dy
\end{eqnarray*}%
which proves the representation formula (\ref{rep0}).

In the previous computations we have used several times Fubini theorem so
we need to prove that some integrability properties hold. Let us suppose
that the support of $f$ is included in $B_{R}(0)$ for some $R>1.$ We denote $%
C_{R}(x)=\{y:\left\vert x\right\vert -R\leq \left\vert y\right\vert \leq
\left\vert x\right\vert +R\}$ and we have $B_{R}(x)\subset C_{R}(x).$ First
of all
$$
\left\vert \phi (x)\partial _{i}Q_{d}(z)\partial
_{i}f(x-z)\right\vert \leq \left\Vert \partial _{i}f\right\Vert _{\infty
}\left\vert \phi (x)\right\vert \left\vert \partial
_{i}Q_{d}(z)1_{C_{R}(x)}(z)\right\vert
$$
and
\begin{equation*}
\int_{C_{R}(x)}\left\vert \partial _{i}Q_{d}(z)\right\vert dz\leq
A_{d}\int_{\left\vert x\right\vert -R}^{\left\vert x\right\vert +R}\frac{r}{%
r^{d}}\times r^{d-1}dr=2RA_{d}.
\end{equation*}%
So%
\begin{equation*}
\int \int \left\vert \phi (x)\partial _{i}Q_{d}(z)\partial
_{i}f(x-z)\right\vert dz\mu (dx)\leq 2RA_{d}\|\partial_if\|_\infty\int \left\vert \phi
(x)\right\vert \mu (dx)<\infty .
\end{equation*}%
Similarly
\begin{eqnarray*}
\int \int \left\vert \partial _{i}^{\mu }\phi (x)\partial
_{i}Q_{d}(z)f(x-z)\right\vert dz\mu (dx)
&=&\int \int \left\vert \partial
_{i}^{\mu }\phi (x)\partial _{i}Q_{d}(x-y)f(y)\right\vert dy\mu (dx) \\
&\leq &2RA_{d}\|f\|_\infty\int \left\vert \partial _{i}^{\mu }\phi (x)\right\vert
\mu(dx)<\infty
\end{eqnarray*}%
so all the needed integrability properties hold and our computation is
correct. In particular we have checked that $\int dyf(y)\int \left\vert
\partial _{i}^{\mu }\phi (x)\partial _{i}Q_{d}(x-y)\right\vert \mu
(dx)<\infty $ for every $f\in C_{c}^{1}(\R^{d})$ so $\int \left\vert \partial
_{i}^{\mu }\phi (x)\partial _{i}Q_{d}(x-y)\right\vert \mu (dx)$ is finite $%
dy $ almost surely.

\smallskip

\textbf{B.}
In order to prove (\ref{rep2}), we write $\partial _{\alpha
}f=\sum_{i=1}^{d}\partial _{i}Q_{d}\ast \partial _{i}\partial _{\alpha }f$.
Now, we use the same chain of equalities as above and we obtain
\begin{eqnarray*}
\int \partial _{\alpha }f(x)p_{\mu _{\phi }}(x)dx
&=&\int \partial _{\alpha
}fd\mu _{\phi }\\
&=&(-1)^{\left\vert \alpha \right\vert }\int
f(y)\Big(-\sum_{i=1}^{d}\int \partial _{i}Q_{d}(x-y)\partial _{(\alpha ,i)}^{\mu
}\phi (x)\mu (dx)\Big)dy
\end{eqnarray*}%
so that $\partial _{\alpha }p_{\mu _{\phi }}(y)=-\sum_{i=1}^{d}\int \partial
_{i}Q_{d}(x-y)\partial _{(\alpha ,i)}^{\mu }\phi (x)\mu (dx).$
As for (\ref{rep3}), we have%
\begin{eqnarray*}
\int \partial _{\alpha }f(x)p_{\mu _{\phi }}(x)dx &=&\int \partial _{\alpha
}f(x)\mu _{\phi }(dx)=\int \partial _{\alpha }f(x)\phi (x)\mu (dx) \\
&=&(-1)^{\left\vert \alpha \right\vert }\int f(x)\partial _{\alpha }^{\mu
}\phi (x)\mu (dx)\\
&=&(-1)^{\left\vert \alpha \right\vert }\int f(x)\partial
_{\alpha }^{\mu }\phi (x)p_{\mu }(x)dx.
\end{eqnarray*}%
$\square $

\begin{remark}
Notice that if $1\in W^{1,1}_\mu$ then for any $f\in C^\infty_c$ one has
$$
\Big|\int\partial_ifd\mu\Big|\leq c_i\|f\|_\infty\quad \mbox{with }
c_i=\|\partial_i^\mu 1\|_{L^1_\mu},
\quad i=1,\ldots ,d.
$$
Now, it is known that the above condition implies the existence
of the density, as proved by Malliavin in \cite{bib:[M]}
(see also D. Nualart \cite{bib:[N]}, Lemma 2.1.1), and
Theorem \ref{th1} gives a new proof including the representation
formula in terms of the Riesz transform.

\end{remark}
\subsection{Estimate of the Riesz transform}

As we will see later on, an important fact is to be able to control
the quantities $\partial _{i}Q_{d},$ and more precisely  $\Theta
_{p}(\mu )$ defined by
\begin{equation}
\Theta _{p}(\mu )=\sup_{a\in \R^{d}}\sum_{i=1}^{d}\Big(\int_{\R^{d}}\left\vert
\partial _{i}Q_{d}(x-a)\right\vert ^{\frac{p}{p-1}}\mu (dx))\Big)^{\frac{p-1%
}{p}}.  \label{ThetaRp}
\end{equation}%
This is the main content of Theorem \ref{th2} below. We begin with two preparatory lemmas.

For a probability measure $\mu $ and probability density $\psi $ (a non
negative measurable function $\psi $ with $\int_{\R^{d}}\psi (x)dx=1)$ we
define the probability measure $\psi \ast \mu $ by%
\begin{equation*}
\int f (x)(\psi \ast \mu )(dx):=\int \psi (x)\int f (x+y)\mu (dy)dx.
\end{equation*}

\begin{lemma}\label{lemma1}
Let $p\geq 1.$ If $1\in W_{\mu }^{1,p}$ then $1\in W_{\psi \ast \mu }^{1,p}$
and $\left\Vert 1\right\Vert _{W_{\psi \ast \mu }^{1,p}}\leq \left\Vert
1\right\Vert _{W_{\mu }^{1,p}}.$
\end{lemma}
\textbf{Proof}. On a probability space $(\Omega ,\mathcal{F},P)$ we consider
two independent random variables $F$ and $\Delta $ such that $F\sim \mu (dx)$
and $\Delta \sim \psi (x)dx.$ Then $F+\Delta \sim (\psi \ast \mu )(dx).$
We define $\theta _{i}(x)=\E(\partial _{i}^{\mu }1(F)\mid F+\Delta =x)$ and we
claim that $\partial _{i}^{\psi \ast \mu }1=\theta _{i}.$ In fact, for $f\in
C_{c}^{1}(\R^{d})$ one has%
$$
\begin{array}{rl}
\displaystyle
-\int \partial _{i}f(x)(\psi \ast \mu )(dx)
=&\displaystyle\!\!\!
-\int dx\psi (x)\int \partial
_{i}f(x+y)\mu (dy) \\
=&\displaystyle\!\!\!
-\int dx\psi (x)\int f(x+y)\partial _{i}^{\mu }1(y)\mu (dy) \\
=&\displaystyle \!\!\!
\E(f(F+\Delta
)\partial _{i}^{\mu }1(F))
=\E(f(F+\Delta )\E(\partial _{i}^{\mu }1(F)\mid F+\Delta ))\\
=&\displaystyle\!\!\!
\E(f(F+\Delta
)\theta _{i}(F+\Delta ))
=\int f(x)\times \theta _{i}(x)(\psi \ast \mu )(dx)
\end{array}%
$$
so $\partial _{i}^{\psi \ast \mu }1=\theta _{i}.$ Moreover%
\begin{eqnarray*}
\int \left\vert \theta _{i}(x)\right\vert ^{p}(\psi \ast \mu )(dx)
&=&\E(\left\vert \theta _{i}(F+\Delta )\right\vert ^{p})=\E(\left\vert
\E(\partial _{i}^{\mu }1(F)\mid F+\Delta )\right\vert ^{p}) \\
&\leq &\E(\left\vert \partial _{i}^{\mu }1(F)\right\vert ^{p})=\int
\left\vert \partial _{i}^{\mu }1(x)\right\vert ^{p}\mu (dx)
\end{eqnarray*}%
so $1\in W_{\psi \ast \mu }^{1,p}$ and $\left\Vert 1\right\Vert _{W_{\psi
\ast \mu }^{1,p}}\leq \left\Vert 1\right\Vert _{W_{\mu }^{1,p}}.$ $\square $

\begin{lemma}\label{lemma2}
Let $p_{n},n\in \N$ be a sequence of probability densities such that $%
\sup_{n}\left\Vert p_{n}\right\Vert _{\infty }$ $=C_{\infty }<\infty .$ Suppose
also that the sequence of probability measures $\mu _{n}(dx)=p_{n}(x)dx,n\in
\N$ converges weakly to a probability measure $\mu .$ Then $\mu (dx)=p(x)dx$
and $\left\Vert p\right\Vert _{\infty }\leq C_{\infty }.$
\end{lemma}

\textbf{Proof}. Since $\int p_{n}^{2}(x)dx\leq C_{\infty }$ the sequence $%
p_{n}$ is bounded in $L^{2}(\R^{d})$ and so it is weakly relative compact.
Passing to a subsequence (which we still denote by $p_{n})$ we may find $%
p\in L^{2}(\R^{d})$ such that $\int p_{n}(x)f(x)dx\rightarrow \int p(x)f(x)dx$
for every $f\in L^{2}(\R^{d}).$ But, if $f\in C_{c}(\R^{d})\subset
L^{2}(\R^{d}) $ then $\int p_{n}(x)f(x)dx\rightarrow \int f(x)\mu (dx)$ so
that $\mu (dx)=p(x)dx.$

Let us now check that $p$ is bounded. Using Mazur's theorem we may construct
a convex combination $q_{n}=\sum_{i=1}^{m_{n}}\lambda _{i}^{n}p_{n+i},$ with
$\lambda _{n}^{i}\geq 0,\sum_{i=1}^{m_{n}}\lambda _{i}^{n}=1,$ such that $%
q_{n}\rightarrow p$ strongly in $L^{2}(\R^{d}).$ Then, passing to a
subsequence, we may assume that $q_{n}\rightarrow p$ almost everywhere. It
follows that $p(x)\leq \sup_{n}q_{n}(x)\leq C_{\infty }$ almost everywhere.
And we may change $p$ on a set of null measure. $\square $

\medskip

We are now able to give our basic estimate of $\Theta _{p}(\mu ).$

\begin{theorem}\label{th2}
Let $p>d$ and let $\mu $ be a probability measure on $\R^{d}$ such that
$1\in W_{\mu }^{1,p}.$ So by Theorem \ref{th1} $\mu (dx)=p_{\mu }(x)dx$ with
\begin{equation}
p_{\mu }(x)=\sum_{i=1}^{d}\int \partial _{i}Q_{d}(y-x)\partial _{i}^{\mu
}1(y)\mu (dy).  \label{rep1}
\end{equation}%
Then%
\begin{equation}\label{est1}
\begin{array}{rl}
i) & \quad \Theta _{p}(\mu ) \leq dK_{d,p}\left\Vert 1\right\Vert
_{W_{\mu }^{1,p}}^{k_{d,p}}\smallskip \\
ii) & \quad \left\Vert p_{\mu }\right\Vert _{\infty } \leq 2dK_{d,p}
\left\Vert 1\right\Vert _{W_{\mu }^{1,p}}^{k_{d,p}+1}.
\end{array}%
\end{equation}
with%
\begin{equation*}
k_{d,p}=\frac{(d-1)p}{p-d}\quad \mbox{and}\quad K_{d,p}=1+2A_{d}^{\frac{p}{%
p-1}}\Big(\frac{p-1}{p-d}\cdot 2d\,A_{d}^{\frac{p}{p-1}}\Big)^{k_{d,p}}.
\end{equation*}
\end{theorem}

\begin{remark}
The inequality in (\ref{est1}) $i)$ gives estimates of the kernels $\partial
_{i}Q_{d},i=1,\ldots ,d$ which appear in the Riesz transform. This is the
crucial point in our approach. In Malliavin and E. Nualart \cite{bib:[M.N]},
 the authors use the Riesz transform
on the sphere and they give estimates of the $L^{p}$ norms of the\
corresponding kernels (which are of course different).
\end{remark}

\textbf{Proof}. We will first prove the theorem under the supplementary
assumption:

\begin{equation*}
\text{(H)}\quad p_{\mu }\text{ is bounded.}
\end{equation*}%
We take $\rho >0$ and notice that if $|x-a|>\rho $ then $%
|\partial _{i}Q_{d}(x-a)|\leq A_{d}\rho ^{-(d-1)}$. Since $p>d$,
for any $a\in \R^{d}$ we have
\begin{equation*}
\begin{array}{rl}
\displaystyle
\int |\partial _{i}Q_{d}(x-a)|^{\frac{p}{p-1}}\mu (dx) & \displaystyle\leq
A_{d}^{\frac{p}{p-1}}\rho ^{-(d-1)\frac{p}{p-1}}+\int_{|x-a|\leq \rho
}|\partial _{i}Q_{d}(x-a)|^{\frac{p}{p-1}}\,p_{\mu }(x)\,dx\smallskip \\
& \displaystyle\leq A_{d}^{\frac{p}{p-1}}\Big[\rho ^{-(d-1)\frac{p}{p-1}%
}+\Vert p_{\mu }\Vert _{\infty }\int_{0}^{\rho }\frac{dr}{r^{\frac{d-1}{p-1}}%
}\Big].%
\end{array}%
\end{equation*}%
This gives
\begin{equation}
\int |\partial _{i}Q_{d}(x-a)|^{\frac{p}{p-1}}\mu (dx)\leq A_{d}^{\frac{p}{%
p-1}}\Big[\rho ^{-(d-1)\frac{p}{p-1}}+\Vert p_{\mu }\Vert _{\infty }\,\frac{%
p-1}{p-d}\,\rho ^{\frac{p-d}{p-1}}\Big]<\infty .  \label{QFeps}
\end{equation}%
We use the representation formula (\ref{rep1}) and H\"{o}lder's inequality
and we obtain%
\begin{eqnarray}
p_{\mu }(x) &=&-\sum_{i=1}^{d}\int \partial _{i}Q_{d}(y-x)\partial _{i}^{\mu
}1(y)\mu (dy) \notag\\
&\leq & \left\Vert 1\right\Vert _{W_{\mu
}^{1,p}}\sum_{i=1}^{d}\Big(\int \left\vert \partial _{i}Q_{d}(y-x)\right\vert ^{%
\frac{p}{p-1}}\mu (dy)\Big)^{\frac{p-1}{p}} \notag \\
&\leq &\left\Vert 1\right\Vert _{W_{\mu }^{1,p}}\Big(d+\sum_{i=1}^{d}\int
\left\vert \partial _{i}Q_{d}(y-x)\right\vert ^{\frac{p}{p-1}}\mu (dy)\Big). \label{var1}
\end{eqnarray}%
By using (\ref{QFeps}), we obtain
\begin{equation*}
\Vert p_{\mu }\Vert _{\infty }\leq d\,\,\left\Vert 1\right\Vert _{W_{\mu
}^{1,p}}\Big(A_{d}^{\frac{p}{p-1}}\Big(\rho ^{-(d-1)\frac{p}{p-1}}+\Vert p_{\mu
}\Vert _{\infty }\,\frac{p-1}{p-d}\,\rho ^{\frac{p-d}{p-1}}\Big)+1\Big)
\end{equation*}%
Choose now $\rho =\rho _{\ast }$, with $\rho _{\ast }$ such that
\begin{equation*}
d\,A_{d}^{\frac{p}{p-1}}\,\left\Vert 1\right\Vert _{W_{\mu }^{1,p}}\frac{p-1%
}{p-d}\,\rho _{\ast }^{\frac{p-d}{p-1}}=\frac{1}{2}
\end{equation*}%
that is,
\begin{equation*}
\rho _{\ast }=\Big(\frac{p-1}{p-d}\cdot 2d\,A_{d}^{\frac{p}{p-1}%
}\,\left\Vert 1\right\Vert _{W_{\mu }^{1,p}}\Big)^{-\frac{p-1}{p-d}}.
\end{equation*}%
Then
\begin{equation*}
\Vert p_{\mu }\Vert _{\infty }\leq 2d\,\,\left\Vert 1\right\Vert _{W_{\mu
}^{1,p}}(A_{d}^{\frac{p}{p-1}}\rho _{\ast }^{-(d-1)\frac{p}{p-1}}+1).
\end{equation*}%
Since $\frac{p-1}{p-d}\,\rho _{\ast }^{\frac{p-d}{p-1}}=(2dA_{d}^{p/(p-1)}%
\left\Vert 1\right\Vert _{W_{\mu }^{1,p}})^{-1}$,
by using (\ref{QFeps}) we obtain%
\begin{equation*}
\begin{array}{rl}
\displaystyle
\int |\partial _{i}Q_{d}(x-a)|^{\frac{p}{p-1}}\mu (dx) & \leq \displaystyle%
1+2A_{d}^{\frac{p}{p-1}}\,\rho _{\ast }^{-(d-1)\frac{p}{p-1}}\smallskip\\
& =\displaystyle1+2A_{d}^{\frac{p}{p-1}}\Big(\frac{p-1}{p-d}\cdot 2d\,A_{d}^{%
\frac{p}{p-1}}\Big)^{\frac{(d-1)p}{p-d}}\cdot \left\Vert 1\right\Vert
_{W_{\mu }^{1,p}}^{\frac{(d-1)p}{p-d}}\smallskip \\
& \leq \displaystyle\Big(1+2A_{d}^{\frac{p}{p-1}}\Big(\frac{p-1}{p-d}\cdot
2d\,A_{d}^{\frac{p}{p-1}}\Big)^{\frac{(d-1)p}{p-d}}\Big)\left\Vert
1\right\Vert _{W_{\mu }^{1,p}}^{k_{d,p}} \\
& =\displaystyle K_{d,p}\cdot \left\Vert 1\right\Vert _{W_{\mu
}^{1,p}}^{k_{d,p}}%
\end{array}%
\end{equation*}%
and finally%
\begin{equation*}
\sum_{i=1}^{d}\int\left\vert \partial _{i}Q_{d}(y-x)
\right\vert ^{\frac{p}{p-1}%
}\mu(dy)\leq d(1+K_{d,p})\left\Vert 1\right\Vert _{W_{\mu
}^{1,p}}^{k_{d,p}}\leq 2dK_{d,p}\left\Vert 1\right\Vert _{W_{\mu
}^{1,p}}^{k_{d,p}}.
\end{equation*}%
Using (\ref{var1}) this gives%
\begin{equation*}
\left\Vert p_{\mu }\right\Vert _{\infty }\leq 2dK_{d,p}\left\Vert
1\right\Vert _{W_{\mu }^{1,p}}^{k_{d,p}+1}.
\end{equation*}

So the theorem is proved under the supplementary assumption (H). We remove
now this assumption. We consider a non negative and continuous function $%
\psi $ such that $\int \psi =1$ and $\psi (x)=0$ for $\left\vert
x\right\vert \geq 1.$ Then we define $\psi _{n}(x)=n^d\psi (nx)$ and $\mu
_{n}=\psi _{n}\ast \mu .$ We have $\mu _{n}(dx)=p_{n}(x)dx$ with $%
p_{n}(x)=\int \psi _{n}(x-y)\mu (dy)\emph{.}$ Using Lemma \ref{lemma1} we have $%
1\in W_{\mu _{n}}^{1,p}$ and $\left\Vert 1\right\Vert _{W_{\mu
_{n}}^{1,p}}\leq \left\Vert 1\right\Vert _{W_{\mu }^{1,p}}<\infty .$ Since $%
p_{n}$ is bounded, $\mu _{n}$ verifies assumption (H) and so, using the
first part of the proof, we obtain
\begin{equation*}
\left\Vert p_{n}\right\Vert _{\infty }\leq 2dK_{d,p}\left\Vert
1\right\Vert _{W_{\mu _{n}}^{1,p}}^{k_{d,p}+1}\leq 2dK_{d,p}
\left\Vert 1\right\Vert _{W_{\mu }^{1,p}}^{k_{d,p}+1}.
\end{equation*}%
Clearly $\mu _{n}\rightarrow \mu $ weakly so, using Lemma \ref{lemma2} we may find
$p$ such that $\mu (dx)=p(x)dx$ and $p$ is bounded. So $\mu $ itself
satisfies (H) and the proof is completed. $\square $

\subsection{Regularity of the density}

Theorem \ref{th1} says that $\mu _{\phi }$ has a
density as soon as $\phi \in W_{\mu }^{1,1}$ - and this does not depend on
the dimension $d$ of the space. But if we want to obtain a continuous or a
derivable density, we need more regularity for $\phi .$ The main instrument
in order to obtain such properties is the classical theorem of Morrey
which we recall now (see  Corollary IX.13  in Brezis  \cite{bib:brezis}).
\begin{theorem}\label{morrey}
Let $u\in W^{1,p}(\R^{d}).$ If $1-\frac{d}{p}>0$ then $u$ is H\"{o}lder
continuous of exponent $q=1-\frac{d}{p}.$ Furthermore suppose that $u\in
W^{m,p}(\R^{d})$ and $m-\frac{d}{p}>0.$ Let $k=[m-\frac{d}{p}]$ be the integer
part of $m-\frac{d}{p}$ and $q=\{m-\frac{d}{p}\}$ the fractional part. If $%
k=0$ then $u$ is H\"{o}lder continuous of exponent $q$ and if $k\geq 1$ then $%
u\in C^{k}$ and for any multi index $\alpha $ with $\left\vert \alpha
\right\vert \leq k$ the derivative $\partial _{\alpha }u$ is H\"{o}lder
continuous of exponent $q:$ for any $x,y\in \R^d$,
$$
|\partial_\alpha u(x)-\partial_\alpha u(y)|\leq C_{d,p}\|u\|_{W^{m,p}(\R^d)}\,|x-y|^q
$$
$C_{d,p}$ being dependent on $d$ and $p$ only.
\end{theorem}
It is clear from Theorem \ref{morrey} that there are two ways to improve the
regularity of $u$: one has to increase $m$\ or/and $p.$ If $\phi \in W_{\mu
}^{m,p}$ for a sufficiently large $m$ then Theorem \ref{th1} already gives us
a differentiable density $p_{\mu _{\phi }}.$ But if we want to keep $m$ low
we have to increase $p.$ And in order to be able to do it the key point is
the estimate for $\Theta
_{p}(\mu )$ given in Theorem \ref{th2}. This is done in  next
Theorem \ref{th3}, where  we
use the following natural notation: we allow a multi index to be equal to
the empty set and for $\alpha=\emptyset$, we set $|\alpha|=0$
and $\partial_\alpha
f:=f$.

\begin{theorem}\label{th3}
We consider some $p>d$ and we suppose that $1\in W_{\mu }^{1,p}.$
For $m\geq 1$, let $\phi \in W_{\mu }^{m,p}$, so that
$\mu _{\phi}(dx)=p_{\mu _{\phi }}(x)dx.$ Then the following
statements hold.

\smallskip

\textbf{A.}
We have
$p_{\mu _{\phi }}\in W^{m,p}$ and
\begin{equation}\label{estf1-ii}
\left\Vert p_{\mu _{\phi }}\right\Vert _{W^{m,p}} \leq
(2dK_{d,p})^{1-1/p}\left\Vert 1\right\Vert _{W_{\mu
}^{1,p}}^{k_{d,p}(1-1/p)}\left\Vert \phi \right\Vert _{W_{\mu }^{m,p}}.
\end{equation}%
Moreover, for any multi index $\alpha$ such that
$0\leq |\alpha|=\ell\leq m-1$,
we have
\begin{equation}\label{estf1-i}
\left\Vert \partial_\alpha p_{\mu _{\phi }}\right\Vert _{\infty } \leq
dK_{d,p}\left\Vert 1\right\Vert _{W_{\mu }^{1,p}}^{k_{d,p}}
\left\Vert \phi
\right\Vert _{W_{\mu }^{\ell+1,p}}.
\end{equation}%

\smallskip

\textbf{B.} We have
 $p_{\mu _{\phi }}\in C^{m-1}$.
Moreover,
for any multi index $\beta$ such that
$0\leq |\beta|=k\leq m-2$, $\partial_\beta p_{\mu_\phi}$ is
Lipschitz continuous: for any $x,y\in \R^d$,
$$
|\partial_\beta p_{\mu_\phi}(x)-\partial_\beta p_{\mu_\phi}(y)|
\leq d^2K_{d,p}\left\Vert 1\right\Vert _{W_{\mu }^{1,p}}^{k_{d,p}}
\left\Vert \phi
\right\Vert _{W_{\mu }^{k+2,p}}\, |x-y|.
$$
And for any multi index $\beta$ such that
$|\beta|= m-1$, $\partial_\beta p_{\mu_\phi}$ is
H\"older continuous of exponent $1-d/p$: for any $x,y\in \R^d$,
$$
\left\vert \partial _{\beta
}p_{\mu _{\phi }}(x)-\partial _{\beta }p_{\mu _{\phi }}(y)\right\vert
\leq
C_{d,p}\left\Vert p_{\mu _{\phi }}\right\Vert _{W^{m,p}}\left\vert x-y\right\vert
^{1-d/p}
$$
$C_{d,p}$ being dependent on $d$ and $p$ only.

\medskip

\textbf{C.} We have
 $W_{\mu }^{m,p}\subset \cap _{\delta >0}W^{m,p}(\{p_{\mu }>\delta
\}).$
\end{theorem}

\textbf{Proof}. \textbf{A. }
We use (\ref{rep3}) (with the notation $\partial_\alpha^\mu \phi
:=\phi$ if $\alpha=\emptyset$) and we obtain%
\begin{equation*}
\int \left\vert \partial _{\alpha }p_{\mu _{\phi }}(x)\right\vert
^{p}dx=\int \left\vert \partial _{\alpha }^{\mu }\phi (x)\right\vert
^{p}\left\vert p_{\mu }(x)\right\vert ^{p}dx\leq \left\Vert p_{\mu
}\right\Vert _{\infty }^{p-1}\int \left\vert \partial _{\alpha }^{\mu }\phi
(x)\right\vert ^{p}p_{\mu }(x)dx.
\end{equation*}%
So $\left\Vert \partial _{\alpha }p_{\mu _{\phi }}\right\Vert _{L^{p}}\leq
\left\Vert p_{\mu }\right\Vert _{\infty }^{1-1/p}\left\Vert \partial
_{\alpha }^{\mu }\phi \right\Vert _{L_{\mu }^{p}}$  and by
using (\ref{est1}) we obtain (\ref{estf1-ii}).
Now, using the representation formula (\ref{rep2}), H\"{o}lder's
inequality and Theorem \ref{th2} we get%
$$
\left\vert \partial_\alpha p_{\mu _{\phi }}(x)\right\vert
\leq
\Theta _{p}(\mu )\sum_{i=1}^d
\|\partial_{(\alpha,i)}^\mu \phi\|_{L^p_\mu}
\leq
dK_{d,p}\left\Vert 1\right\Vert
_{W_{\mu }^{1,p}}^{k_{d,p}}
\|\phi\|_{W^{\ell+1,p}_\mu}
$$
and (\ref{estf1-i}) is proved.

\medskip

\textbf{B.}
The fact that $p_{\mu _{\phi }}\in
C^{m-1}(\R^{d})$ and the H\"{o}lder property are standard consequences
of $%
p_{\mu _{\phi }}\in W^{m,p}$, as stated in Theorem \ref{morrey}.
As for the Lipschitz property, it immediately
follows from (\ref{estf1-i}) and
the fact that if $f\in C^1$ with $\|\nabla f\|_\infty
<\infty$ then $|f(x)-f(y)|\leq \sum_{i=1}^d \|\partial_i f\|_\infty\,
|x-y|$.

\medskip

\textbf{C. }
We have $p_{\mu _{\phi }}(x)dx=\phi (x)\mu (dx)=\phi (x)p_{\mu }(x)dx$ so
$\phi (x)=p_{\mu _{\phi }}(x)/p_{\mu }(x)$ if $p_{\mu }(x)>0.$ And since $%
p_{\mu _{\phi }},p_{\mu }\in W^{m,p}(\R^{d})$ we obtain $\phi \in
W^{m,p}(\{p_{\mu }>\delta \}).$ $\square $

\subsection{Estimate of the tails of the density}

In order to study the behavior of the tails of the density, we need the following computational rules.

\begin{lemma}\label{lemma-rules}
\textbf{A. }
If $\phi \in W_{\mu }^{1,p}$ and $\psi \in
C_{b}^{1}(\R^{d})$ then $\psi \phi \in W_{\mu }^{1,p}$ and
\begin{equation}
\partial _{i}^{\mu }(\psi \phi )=\psi \partial _{i}^{\mu }\phi +\partial
_{i}\psi \phi .  \label{dist3}
\end{equation}%
In particular, if $1 \in W_{\mu }^{1,p}$ then for any $\psi \in
C_{b}^{1}(\R^{d})$ one has  $\psi \in W_{\mu }^{1,p}$ and
\begin{equation}
\partial _{i}^{\mu }\psi =\psi \partial _{i}^{\mu }1 +\partial
_{i}\psi.  \label{dist3-bis}
\end{equation}%
\textbf{B. }
Suppose that $1\in W_{\mu }^{1,p}$. If $\psi \in C_{b}^{1}(\R^{m})$ and if
$u=(u_{1},\ldots  ,u_{m})\,:\,\R^{d}\rightarrow \R^{m}$ is such that $u_{j}\in
C_{b}^{1}(\R^{d})$, $j=1,\ldots  ,m$, then $\psi \circ u\in W_{\mu }^{1,p}$
and
\begin{equation*}
\partial _{i}^{\mu }\psi \circ u=\sum_{j=1}^{m}(\partial _{j}\psi )\circ
u\partial _{i}^{\mu }u_{j}+T_{\psi }\circ u\,\partial _{i}^{\mu }1
\end{equation*}%
where
\begin{equation*}
T_{\psi }(x)=\psi (x)-\sum_{j=1}^{d}\partial _{j}\psi (x)x_{j}.
\end{equation*}
\end{lemma}

\textbf{Proof}. \textbf{A. }
Since $\psi $ and $\partial _{i}\psi $ are bounded,
$\psi
\phi,\psi \partial _{i}^{\mu }\phi ,\partial _{i}\psi
\phi \in L_{\mu }^{p}.$ So we just have to check the integration by parts
formula. We have%
$$
\int \partial _{i}f\psi \phi \,d\mu =\int \partial _{i}(f\psi )\phi \,d\mu
-\int f\partial _{i}\psi \phi \,d\mu
=-\int f\psi \partial _{i}^{\mu }\phi \,d\mu -\int f\partial _{i}\psi \phi
\,d\mu
$$
and the statement holds.

\smallskip
\textbf{B. }
By using (\ref{dist3-bis}), we have
\begin{equation*}
\partial _{i}^{\mu }(\psi \circ u)=\psi \circ u\,\partial _{i}^{\mu
}1+\partial _{i}(\psi \circ u)=\psi \circ u\,\partial _{i}^{\mu
}1+\sum_{j=1}^{m}(\partial _{j}\psi )\circ u\,\partial _{i}u_{j}.
\end{equation*}%
The formula now follows by
inserting $\partial _{i}u_{j}=\partial _{i}^{\mu
}u_{j}-u_{j}\partial _{i}^{\mu }1,$ as given by (\ref{dist3-bis}). $\square $

\medskip

We give now a result which allows to estimate the queues of $p_{\mu}.$

\begin{proposition}\label{queue}
Let $\phi \in C_{b}^{1}(\R^{d})$ be a function such that $%
1_{B_{1}(0)}\leq \phi \leq 1_{B_{2}(0)}$ and $|\nabla \phi|\leq 1$.
We set  $\phi _{x}(y)=\phi
(x-y)$ and we assume that $1\in W_{\mu }^{1,p}$ with $p>d,$ so, in view of
Lemma \ref{lemma-rules}, $\phi _{x}\in W_{\mu }^{1,p}$. Then we have the
representation%
\begin{equation*}
p_{\mu }(x)=\sum_{i=1}^{d}\int \partial _{i}Q_{d}(y-x)\partial _{i}^{\mu
}\phi _{x}(y)1_{\{\left\vert y-x\right\vert <2\}}\mu (dy).
\end{equation*}%
As a consequence, for any positive $a<\frac{1}{d}-\frac{1}{p}$ one has
\begin{equation}\label{est-queue}
p_{\mu }(x)\leq \Theta _{\overline{p}}(\mu )\big(d+\|1\|_{W^{1,p}_\mu}\big)
\,\mu \big(B_{2}(x)\big)^{a}.
\end{equation}%
where $\overline{p}=1/(a+\frac{1}{p})$. In particular,
\begin{equation}\label{est-queue2}
\lim_{|x|\to\infty} p_\mu(x)=0.
\end{equation}
\end{proposition}

\textbf{Proof}. By Lemma \ref{lemma-rules}, $\phi _{x}\in W_{\mu }^{1,p}$
and $\partial _{i}^{\mu }\phi _{x}=\phi _{x}\partial _{i}^{\mu
}1+\partial _{i}\phi _{x}$, so that $\partial _{i}^{\mu }\phi _{x}=\partial
_{i}^{\mu }\phi _{x}1_{B_{2}(x)}$. Now, for $f\in C_{c}^{1}(B_{1}(x))$ we
have
\begin{align*}
\int f(y)\mu(dy)& =\int f(y) \phi _{x}(y)\mu(dy)=\int f(y)p_{\mu _{\phi _{x}}}(y)dy \\
& =-\int f(y)\sum_{i=1}^{d}\int \partial _{i}Q_{d}(z-y)\partial _{i}^{\mu
}\phi _{x}(z)\mu (dz)dy\\
& =-\int f(y)\sum_{i=1}^{d}\int \partial _{i}Q_{d}(z-y)\partial _{i}^{\mu
}\phi _{x}(z)1_{B_{2}(x)}(z)\mu (dz)dy.
\end{align*}%
It follows that for $y\in B_{1}(x)$ we have
$$
p_{\mu }(y)
=-\sum_{i=1}^{d}\int \partial _{i}Q_{d}(z-y)\partial _{i}^{\mu }\phi
_{x}(z)1_{B_{2}(x)}(z)\mu (dz).
$$
We consider now $y=x$ and we take $a\in (0,\frac{1}{d}-\frac{1}{p}).$ Using H%
\"{o}lder's inequality we obtain%
$$
p_{\mu }(x) \leq \mu (B_{2}(x))^{a}\sum_{i=1}^{d}I_{i}\quad \mbox{with }
I_{i} =\Big(\int \left\vert \partial _{i}Q_{d}(z-x)\partial _{i}^{\mu }\phi
_{x}(z)\right\vert ^{\frac{1}{1-a}}\mu (dz)\Big)^{a-1}.
$$
Notice that $1<d(1-a)/(1-da)<p(1-a).$ We take $\beta $ such that $%
d(1-a)/(1-da)<\beta <p(1-a)$ and we denote by $\alpha $ the conjugate of $%
\beta .$ Using again H\"{o}lder's inequality we obtain%
\begin{eqnarray*}
I_{i} &\leq &\Big(\int \left\vert \partial _{i}Q_{d}(z-x)\right\vert ^{\frac{%
\alpha }{1-a}}\mu (dz)\Big)^{(a-1)/\alpha }\Big(\int \left\vert \partial _{i}^{\mu
}\phi _{x}(z)\right\vert ^{\frac{\beta }{1-a}}\mu (dz)\Big)^{(a-1)/\beta } \\
&\leq &\Big(\int \left\vert \partial _{i}Q_{d}(z-x)\right\vert ^{\frac{\alpha }{%
1-a}}\mu (dz)\Big)^{(a-1)/\alpha }\left\Vert \partial _{i}^{\mu }\phi
_{x}\right\Vert _{L_{\mu }^{p}}.
\end{eqnarray*}%
We let $\beta \uparrow p(1-a)$ so that
\begin{equation*}
\frac{\alpha }{1-a}=\frac{\beta }{(\beta -1)(1-a)}\rightarrow \frac{p}{%
p(1-a)-1}=\frac{\overline{p}}{\overline{p}-1}.
\end{equation*}%
So we obtain
\begin{equation*}
I_{i}\leq \Theta _{\overline{p}}(\mu )\left\Vert \partial _{i}^{\mu }\phi
_{x}\right\Vert _{L_{\mu }^{p}}
\end{equation*}%
and then
$$
p_{\mu }(x)
\leq
\Theta _{\overline{p}}(\mu )\left\Vert \phi
_{x}\right\Vert _{W^{1,p}_{\mu }}
\mu (B_{2}(x))^{a}.
$$
Now, since
$\partial _{i}^{\mu }\phi _{x}=\phi _{x}\partial _{i}^{\mu
}1+\partial _{i}\phi _{x}$, we have $\left\Vert \phi
_{x}\right\Vert _{W^{1,p}_{\mu }}\leq d+\|1\|_{W^{1,p}_\mu}$
and (\ref{est-queue}) is proved.
Finally, $\I_{B_{2}(x)}\to 0$ a.s. when $|x|\to\infty$ and by using the Lebesgue dominated convergence
theorem, one has
$\mu \big(B_{2}(x)\big)=\int\I_{B_{2}(x)}(y)\mu(dy)\to 0$.
By applying (\ref{est-queue}), one obtains (\ref{est-queue2}).
$\square $

\subsection{On the set of strict positivity for the density}

Suppose that $1\in W^{2,1}_\mu$ and set $U_\mu=\{p_\mu>0\}$. We define the matrix field
$C_\mu\,:\,\R^d\to \R^d\times \R^d$ through
$$
C_\mu^{ij}(x)=\partial_i^\mu1(x)\partial_j^\mu1(x),\quad i,j=1,\ldots,d.
$$
For $x,y\in \R^d$, we set
\begin{eqnarray}
A^\mu_{x,y}&=&\displaystyle
\Big\{\varphi\in C^1([0,1],U_\mu)\,;\,\varphi_0=x,\varphi_1=y\Big\}\nonumber\\
d_\mu(x,y)
&=&\displaystyle
\inf\Big\{\int_0^1\langle C_\mu(\varphi_t)\dot \varphi_t,\dot \varphi_t\rangle^{1/2} dt\,;\,
\varphi\in A^\mu_{x,y}\Big\}\label{dist-1}
\end{eqnarray}
with the understanding $d_\mu(x,y)=+\infty$
if $A_{x,y}=\emptyset$. Notice that
$d_\mu(x,y)=+\infty$ if $x$ and $y$ belong to two different
connected components of the open set $U_\mu$, if they exist. Moreover, it is easy to see that
$d_\mu(x,y)$ does not define in general a distance but
only a semi-distance. In fact, as an example, take $p_\mu$ as a smooth probability density on $\R$
which is constant on some interval $(a,b)$, $a<b$. Then, $\partial^\mu 1=\partial\ln p_\mu
\equiv 0$ on $(a,b)$, so that $d_\mu(x,y)=0$ for any $x,y\in (a,b)$.

Then we have the following representation formula and estimates for the density.

\begin{proposition}\label{bell-prop}
Let $1\in W^{2,p}_\mu$ with $p>d$, and let $x,x_0\in U_\mu$ be such that $A^\mu_{x_0,x}\neq \emptyset.$
Then, for any $\varphi\in A^\mu_{x_0,x}$ one has
\begin{equation}\label{bell-formula}
p_\mu(x)=p_\mu(x_0)\exp\Big(\int_0^1\langle\partial^\mu1(\varphi_t), \dot\varphi_t\rangle\,dt\Big)
\end{equation}
As a consequence,
\begin{equation}\label{dist-3}
p_\mu(x_0)e^{-d_\mu(x_0,x)}\leq
p_\mu(x) \leq
p_\mu(x_0)e^{d_\mu(x_0,x)}
.
\end{equation}

\end{proposition}

\textbf{Proof.}
If $1\in W^{1,p}_\mu$ with $p>d$, the density $p_\mu$ exists and is continuous, so that
$U_\mu= \{p_\mu>0\}$ is open. Now, for $\varphi\in A^\mu_{x_0,x}$,
 one has
$\partial_i\ln p_\mu (\varphi_{x_0,x}(t))=\partial_i^\mu1(\varphi_{x_0,x}(t))$
so that
$$
d\ln p_\mu (\varphi(t))
=\sum_{i=1}^d \partial_i^\mu1(\varphi(t))\dot\varphi^i(t)dt
=\langle\partial^\mu1(\varphi_t), \dot\varphi_t\rangle
$$
and (\ref{bell-formula}) follows by integrating over $[0,1]$. Now, for any $\varphi\in A^\mu_{x_0,x}$
one has
$$
\Big|\ln \frac{p(x)}{p(x_0)}\Big|
=\Big|\int_0^1\langle\partial^\mu1(\varphi_t), \dot\varphi_t\rangle\,dt\Big|
\leq
\int_0^1|\langle\partial^\mu1(\varphi_t), \dot\varphi_t\rangle|
\,dt
=\int_0^1\langle C_\mu(\varphi_t)\dot\varphi_t, \dot\varphi_t\rangle^{1/2}
\,dt.
$$
By taking the $\inf$ over $A^\mu_{x_0,x}$ one proves (\ref{dist-3}).
$\square$

\medskip

We can now state the main result of this section.

\begin{proposition}\label{mall-conj-prop}
Suppose that $1\in W^{2,p}_\mu$ with $p>d$. Then the following statements hold.

\medskip

$i)$ If $p_\mu(x_n)\to 0$ then $d_\mu(x_n,x_1) \to\infty$.

\medskip

$ii)$ If $U_\mu$ is connected, the converse of $i)$ holds:
if $d_\mu(x_n,x_1) \to\infty$ then $p_\mu(x_n)\to 0$.

\medskip

$iii)$
$\partial^\mu 1$ is locally
bounded (that is, bounded on compact sets of $\R^d$)
if and only if $\{p_\mu>0\}=\R^d$.
\end{proposition}

\textbf{Proof.}
$i)$ The statement immediately follows by the first inequality in (\ref{dist-3}).

\smallskip

$ii)$ By contradiction, we assume that $p_\mu(x_n)\nrightarrow 0$: there exist $ c>0$
and a subsequence
$\{x_{n_k}\}_k$ such that $p_\mu(x_{n_k})\geq c$ for any $k$. By
Proposition \ref{queue}, and in particular  (\ref{est-queue2}),
there exists $R>0$ such that $p_\mu(x)<c$ if $|x|>R$. This gives that
the sequence $\{x_{n_k}\}_k$ is bounded and then there exists a
further subsequence $\{x_{n_{k_\ell}}\}_\ell$ converging to some point $\bar x$.
Now, since $p_\mu$ is continuous and $p_\mu(\bar x)\geq c$, there exists $r>0$ such that
$p_\mu\geq \frac c 2$ in the ball $B(\bar x,r)$, which of course contains the points
$x_{n_{k_\ell}}$ for any large $\ell$. This means that the path $\varphi^\ell$
joining $\bar x$ to $x_{n_{k_\ell}}$ at constant speed belong to $A_{\bar x, x_{n_{k_\ell}}}$
for any large $\ell$. Therefore,
$$
\int_0^1|\langle\partial^\mu1(\varphi^\ell_t), \dot\varphi^\ell_t\rangle|
\,dt
\leq
\int_0^1|\partial^\mu1(\varphi^\ell_t)|| \dot\varphi^\ell_t|
\,dt
\leq C\times |\bar x- x_{n_{k_\ell}}|
\leq C\times r
$$
in which we have used the fact that $\partial^\mu1$ is bounded on $B(\bar x,r)$.
It follows that for some $\ell_0$,
$$
\sup_{\ell\geq \ell_0}d_\mu(\bar x, x_{n_{k_\ell}})\leq  C\times r.
$$
Now,
$$
\sup_{\ell\geq \ell_0}d_\mu(x_1, x_{n_{k_\ell}})\leq  d_\mu(x_1, \bar x)+
C\times r
$$
and $d_\mu(x_1, \bar x)<\infty$ because $U_\mu$ is connected, and this gives a contradiction.

\smallskip

$iii)$
If $\partial^\mu 1$ is bounded on compact sets of $\R^d$, then for $x\in U_\mu$
by (\ref{bell-formula}) we get
$$
p_\mu(x)\leq p_\mu(x_0)\exp\Big(C\sup_{t\in[0,1]}|\dot\varphi_{x_0,x}|\Big)
$$
where $x_0\in U_\mu$ is such that $A^\mu_{x_0,x}\neq \emptyset$.
Now, if $\partial U_\mu\neq \emptyset$, we can
let $x_0$ tend to the boundary of $U_\mu$ and in such a case we obtain
$p_\mu(x)=0$, which is a contradiction. Therefore, $p_\mu>0$ everywhere.
On the contrary, it is sufficient to recall that
$\partial^\mu1(x)=\partial\ln p_\mu$
is continuous on $U_\mu=\R^d$.
$\square$

\begin{remark}\label{bell-rem}
Parts $i)$ and $ii)$ of above Proposition \ref{mall-conj-prop} allow to discuss
the Malliavin conjecture about the set of the strict positivity
points of the density, as already described in the Introduction at page \pageref{mall-conj}.
For further details, we address to next Proposition \ref{prop-conj}.

Furthermore,  $iii)$ says that if $1\in W^{1,p}_\mu$ with $\partial_i^\mu 1$ locally bounded
then we can take $x_0=0$ and $\varphi_{x_0,x}(t)=tx$, so that
$$
p_\mu(x)=p_\mu(0)\exp\Big(\int_0^1\sum_{i=1}^d x_i\partial_i^\mu1(tx)\, dt\Big).
$$
Such a representation formula has been already given by Bell in \cite{bib:bell}.
\end{remark}

\begin{remark}\label{bell-rem-energy}
It is easy to see that all the results of this section hold
if the semi-distance $d_f$ is replaced
by the square root of the energy associated to the matrix field $C_\mu$, which is defined by
\begin{equation}\label{dist-2bis}
\overline{d}_\mu(x,y)
=
\inf\Big\{\Big(\int_0^1\langle C_\mu(\varphi_t)\dot \varphi_t,\dot \varphi_t\rangle dt\Big)^{1/2}\,;\,
\varphi\in A_{x,y}\Big\}
\end{equation}
with $\overline{d}_\mu(x,y)=+\infty$
if $A_{x,y}=\emptyset$. Again, $\overline{d}_\mu(x,y)$
defines only a semi-distance and one has $d_\mu(x,y)\leq \overline{d}_\mu(x,y)$.
\end{remark}

\subsection{Local integration by parts formulas}

The assumptions in the previous sections are global - and this may fail in
many interesting cases - for example for diffusion processes living in a
region of the space or, as a more elementary example, for the exponential
distribution. So in this section we give a hint about the localized version
of the results presented above.

An open domain $D\subset R^{d}$ is given. We
recall that $L_{\mu }(D)=\{f:\int_{D}\vert f(x)\vert ^{p}d\mu
(x)<\infty \}$ and $W_{\mu }^{1,p}(D)$ is the space of the functions $\phi
\in L_{\mu }(D)$ which verify the integration by parts formula $\int \phi
\partial _{i}fd\mu =-\int \theta _{i}fd\mu $ for test functions $f$ which
have a compact support included in $D.$ And $\partial^{\mu,D} _{i}\phi :=\theta
_{i}\in L_{\mu }(D).$ The space $W_{\mu }^{m,p}(D)$ is similarly defined.
Our aim is to give sufficient conditions in order that $\mu $
has a smooth density on $D,$ that means that we look for a smooth function $%
p $ such that $\int_{D}f(x)d\mu (x)=\int_{D}f(x)p(x)dx.$ And we want to give
estimates for $p$ and its derivatives in terms of the Sobolev norms of $%
W_{\mu }^{m,p}(D).$ 

The main step in our approach is a truncation argument that we present now.
Given $-\infty \leq a\leq b\leq \infty $ and $\varepsilon >0$\ we define $\psi
_{\varepsilon ,a,b}:\R\rightarrow \R_{+}$ by
\[
\psi _{\varepsilon ,a,b}(x)=1_{(a-\varepsilon
,a]}(x)\exp \Big(1-\frac{\varepsilon }{x+\varepsilon -a}\Big) +1_{(a,b)}(x)+
1_{[b,b+\varepsilon )}(x)\exp \Big(1-\frac{\varepsilon }{x-b-\varepsilon }\Big)
\]%
with the convention $1_{(a-\varepsilon ,a)}=0$ if $a=-\infty $ and $%
1_{(b,b+\varepsilon )}=0$ if $b=\infty .$ Notice that $\psi _{\varepsilon
,a,b}\in C_{b}^{\infty }(\R)$ and
\[
\sup_{x\in (a-\varepsilon,b+\varepsilon)}\left\vert \partial _{x}\ln \psi _{\varepsilon
,a,b}(x)\right\vert ^{p}\psi _{\varepsilon ,a,b}(x)\leq \varepsilon
^{-p}\sup_{y>0}y^{2p}e^{-y}.
\]%
For $x=(x^{1},\ldots ,x^{d})$ and $i\in \{1,\ldots ,d\}$ we denote $\widehat{x}%
_{i}=(x^{1},\ldots ,x^{i-1},x^{i+1},\ldots ,x^{d})$ and for $y\in \R$ we put $(%
\widehat{x}_{i},y)=(x^{1},\ldots ,x^{i-1},y,x^{i+1},\ldots ,x^{d}).$ Then we define
\[
a(\widehat{x}_{i})=\inf_{y\in \R}\Big\{y:d((\widehat{x}_{i},y),D^{c})>2%
\varepsilon \Big\},
\quad
b(\widehat{x}_{i})=\sup_{y\in \R}\Big\{y:d((\widehat{x}%
_{i},y),D^{c})>2\varepsilon \Big\}
\]%
with the convention $a(\widehat{x}_{i})=b(\widehat{x}_{i})=0$ if $\{y:d((%
\widehat{x}_{i},y),D^{c})>2\varepsilon \}=\emptyset .$ Finally we define%
\[
\Psi _{D,\varepsilon }(x)=\prod_{i=1}^{d}\psi _{\varepsilon ,a(\widehat{x}%
_{i}),b(\widehat{x}_{i})}(x_{i}).
\]%
We denote $D_{\varepsilon }=\{x:d(x,D^{c})\geq \varepsilon \}$ so that $%
1_{D_{2\varepsilon }}\leq \Psi _{D,\varepsilon }\leq 1_{D_{\varepsilon }}.$
And we also have
\begin{equation}
\sup_{x\in D_\varepsilon}\left\vert \partial _{x}\ln \Psi _{D,\varepsilon
}(x)\right\vert ^{p}\Psi _{D,\varepsilon }(x)\leq d\varepsilon
^{-p}\sup_{y>0}y^{2p}e^{-y}.  \label{LOC1}
\end{equation}

We are now able to give the main result in this section. The symbol
$\nu|_{D}$ denotes the measure $\nu$ restricted to the open set $D$.

\begin{theorem}
\textbf{A}. Suppose that $\phi \in W_{\mu }^{1,1}(D).$ Then $\mu _{\phi}|_{D}
(dx)=p_{\mu _{\phi }}(x)dx.$

\textbf{B}. Suppose that $1\in W_{\mu }^{1,p}(D)$ for some $p>d.$ Then for
each $\varepsilon >0$%
\[
\sup_{x\in \R^{d}}\sum_{i=1}^{d}\int_{D_{2\varepsilon }}\vert \partial
_{i}Q_{d}(y-x)\vert ^{p/(p-1)}\mu (dy)\leq C_{d,p}\varepsilon
^{-p}\Vert 1\Vert _{W_{\mu }^{1,p}(D)}.
\]

\textbf{C}. Suppose that $1\in W_{\mu }^{1,p}(D)$ for some $p>d.$ Then for $%
\phi \in W_{\mu }^{m,p}(D)$ we have $\mu _{\phi }|_{D}(dx)=p_{\phi, D
}(x)dx$ and%
\[
p_{\phi, D }(x)=-\sum_{i=1}^{d}\int \partial _{i}Q_{d}(y-x)(\Psi
_{D,\varepsilon }\partial _{i}^{\mu }\phi +\phi \partial _{i}\Psi
_{D,\varepsilon })\mu (dx)\quad \mbox{for}\quad x\in D_{\varepsilon }.
\]%
Moreover $p_{\phi,D }\in \cap _{\varepsilon >0}W^{m,p}(D_{\varepsilon })
$ and
\[
\Vert p_{\phi,D }\Vert _{W^{m,p}(D_{\varepsilon })}\leq
C_{p,d}(1+\varepsilon ^{-1})(1\vee \Vert 1\Vert _{W_{\mu
}^{1,p}(D)})^{k_{d,p}(1-1/p)}\Vert \phi \Vert _{W_{\mu
}^{m,p}(D)}.
\]%
Finally, $p_{\phi,D }$ is $m-1$ times differentiable on $D$ and
for every multi-index $\alpha $ of length $0\leq \ell\leq m-1$ one has
\[
\Vert \partial _{\alpha }p_{\phi,D }\Vert _{\infty }\leq
C_{p,d}(1+\varepsilon ^{-\ell})(1\vee \Vert 1\Vert _{W_{\mu
}^{1,p}(D)})^{k_{d,p}(1-1/p)}\Vert \phi \Vert _{W_{\mu
}^{l+1,p}(D)}.
\]

\end{theorem}

\textbf{Proof}.
We denote $\mu _{D,\varepsilon }(dx)=\Psi _{D,\varepsilon
}(x)\mu (dx).$ Let us first show that if $\phi \in W_{\mu }^{1,p}(D)$
then $\phi \in W_{\mu_{D,\varepsilon} }^{1,p}(\R^d).$
In fact, if $f\in C_c^\infty(\R^d)$ then
$f\Psi _{D,\varepsilon}\in C_c^\infty(D)$ and similarly to what developed in
Lemma \ref{lemma-rules}, one has
\begin{eqnarray*}
\int \partial _{i}f\phi d\mu _{D,\varepsilon } &=&\int \partial _{i}f\phi\Psi _{D,\varepsilon
}d\mu =-\int \partial _{i}^{\mu }(\Psi _{D,\varepsilon
}\phi )fd\mu  \\
&=&-\int (\Psi _{D,\varepsilon }\partial _{i}^{\mu }\phi +\phi
\partial _{i}\Psi _{D,\varepsilon })fd\mu
=-\int (\partial _{i}^{\mu }\phi +\phi \partial _{i}\ln \Psi
_{D,\varepsilon })fd\mu _{D,\varepsilon }
\end{eqnarray*}%
so that $\partial _{i}^{\mu _{D,\varepsilon }}\phi =\partial _{i}^{\mu }\phi
+\phi \partial _{i}\ln \Psi _{D,\varepsilon }.$ And by using (\ref{LOC1}) we have $%
\partial _{i}^{\mu _{D,\varepsilon }}\phi \in L_{\mu _{D,\varepsilon
}}^{p}(\R^{d})$:
\begin{eqnarray*}
\sum_{i=1}^{d}\int \left\vert \partial _{i}^{\mu _{D,\varepsilon }}\phi
\right\vert ^{p}d\mu _{D,\varepsilon } &=&\sum_{i=1}^{d}\int \left\vert
\partial _{i}^{\mu }\phi +\phi \partial _{i}\ln \Psi _{D,\varepsilon
}\right\vert ^{p}\Psi _{D,\varepsilon }d\mu  \\
&\leq &\sum_{i=1}^{d}\int_{D}\left\vert \partial _{i}^{\mu }\phi +\phi
\partial _{i}\ln \Psi _{D,\varepsilon }\right\vert ^{p}\Psi _{D,\varepsilon
}d\mu \leq C_{p}\varepsilon ^{-p}\left\Vert \phi \right\Vert _{W_{\mu
}^{1,p}(D)}^{p}.
\end{eqnarray*}%
It follows that
\[
\left\Vert \phi \right\Vert _{W_{\mu _{D,\varepsilon }}^{1,p}(R^{d})}\leq
C_{p}\varepsilon ^{-1}\left\Vert \phi \right\Vert _{W_{\mu }^{1,p}(D)}.
\]
Setting $\mu _{D,\varepsilon ,\phi }(dx):=\phi \mu _{D,\varepsilon }(dx)
$, we can use Theorem \ref{th3} and we obtain $\mu _{D,\varepsilon ,\phi
}(dx)=p_{D,\varepsilon ,\phi }(x)dx$ with $p_{D,\varepsilon ,\phi }\in
W^{1,p}(\R^{d}).$ Similarly we prove that if $\phi \in W_{\mu }^{m,p}(D)$ then $p_{D,\varepsilon
,\phi }\in W^{m,p}(\R^{d}).$ We notice that for a function $f$ with the
support included in $D_{2\varepsilon }$ we have $\int f\phi d\mu =\int f\phi d\mu
_{D,\varepsilon }.$ It follows that $\mu_\phi| _{D_{2\varepsilon
}}(dx)=p_{D,\varepsilon }(x)dx.$

\smallskip

Now, statement \textbf{A.} immediately follows from the above arguments.

\smallskip

\textbf{B.}
Suppose now that $1\in W_{\mu }^{1,p}(D)$ for some $p>d.$ Then%
\begin{eqnarray*}
\int_{D_{2\varepsilon }}\left\vert \partial _{i}Q_{d}(x-y)\right\vert
^{p/(p-1)}\mu (dy) &\leq &\int \left\vert \partial _{i}Q_{d}(x-y)\right\vert
^{p/(p-1)}\Psi _{D,\varepsilon }(y)\mu (dy) \\
&=&\int \left\vert \partial _{i}Q_{d}(x-y)\right\vert ^{p/(p-1)}\mu
_{D,\varepsilon }(dy) \\
&\leq &dK_{d,p}(1\vee \left\Vert 1\right\Vert _{W_{\mu _{D,\varepsilon
}}^{1,p}})^{k_{d,p}}\leq C_{d,p}\varepsilon ^{-p}\left\Vert 1\right\Vert
_{W_{\mu }^{1,p}(D)}.
\end{eqnarray*}%
The upper bounds for the density and its derivatives are proved in a similar
way.

\smallskip

Finally, \textbf{C} follows similarly as in $iii)$ of Proposition \ref{mall-conj-prop}.
$\square $

\section{Integration by parts formulas for random variables and Riesz
transform}\label{sect-IBP}

Let $(\Omega ,\mathcal{F},P)$ ba a probability space and let $F$ and $G$
be two random variables taking values in $\R^d$ and $\R$ respectively.

\begin{definition}
Given a multi-index $\alpha $ and a power $p\geq 1$, we say that the
integration by parts formula $IP_{\alpha ,p}(F,G)$ holds if there exists a
random variable $H_{\alpha }(F;G)\in L^{p}$ such that%
\begin{equation}
IP_{\alpha ,p}(F,G)\quad \E(\partial _{\alpha }f(F)G)=(-1)^{\left\vert \alpha
\right\vert }\E(f(F)H_{\alpha }(F;G)),\quad \forall f\in C_{c}^{\left\vert
\alpha \right\vert }(\R^{d}).  \label{den1}
\end{equation}%
We define $W_{F}^{m,p}$ to be the space of the random variables $G\in
L^{p}$ such that $IP_{\alpha ,p}(F,G)$ holds for every multi index $\alpha $
with $\left\vert \alpha \right\vert \leq m.$ For $G\in W_{F}^{m,p}$ we
define
$$
\partial _{\alpha }^{F}G=\E(H_{\alpha }(F;G)\mid F).
$$
\end{definition}

We denote by $\mu _{F}$ the law of $F$ and $\mu _{F,G}(f):=\E(f(F)G).$ So $%
\mu _{F}=\mu _{F,1}.$ Moreover we denote $\phi _{G}(x)=\E(G\mid F=x).$ Then
it is easy to check that

\begin{equation*}
W_{F}^{m,p}=\{G\in L^{p}:\phi _{G}\in W_{\mu _{F}}^{m,p}\}\quad and\quad
\partial _{\alpha }^{F}G=\partial _{\alpha }^{\mu _{F}}\phi _{G}(F).
\end{equation*}%
We also define the norms
\begin{equation*}
\left\Vert G\right\Vert _{W_{F}^{m,p}}=\left\Vert G\right\Vert
_{L^{p}}+\sum_{1\leq \left\vert \alpha \right\vert \leq m}(\E(\left\vert
\partial _{\alpha }^{F}G\right\vert ^{p}))^{1/p}.
\end{equation*}%
It is easy to see that $(W_{F}^{m,p},\left\Vert \cdot \right\Vert
_{W_{F}^{m,p}})$ is a Banach space.

\begin{remark}
Notice that $\E(\left\vert \partial _{\alpha }^{F}G\right\vert ^{p})\leq
\E(\left\vert H_{\alpha }(F;G)\right\vert ^{p})$ so that $\partial _{\alpha
}^{F}G$ is the weight of minimal $L^{p}$ norm which verifies $IP_{\alpha
,p}(F;G).$ In particular%
\begin{equation*}
\left\Vert G\right\Vert _{W_{F}^{m,p}}\leq \left\Vert G\right\Vert
_{L^{p}}+\sum_{1\leq \left\vert \alpha \right\vert \leq m}\left\Vert
H_{\alpha }(F;G)\right\Vert _{L^{p}}
\end{equation*}%
and this last quantity is the one which naturally appears in concrete
computations.
\end{remark}

We can resume the result of Section \ref{sect-sobolev} as follows. As for the
density, we obtain

\begin{theorem}\label{th4}
\textbf{A. }
Suppose that $1\in W_{F}^{1,p}$ for some $p>d.$ Then $\mu
_{F}(dx)=p_{F}(x)dx$ and $p_{F}\in C_{b}(\R^{d}).$ Moreover%
\begin{eqnarray*}
\Theta _{F}(p) &:=&\sup_{a\in \R^{d}}\sum_{i=1}^{d}\Big(\E(\left\vert \partial
_{i}Q_{d}(F-a)\right\vert ^{p/(p-1)})\Big)^{(p-1)/p}\leq K_{d,p}\left\Vert 1\right\Vert
_{W_{F}^{1,p}}^{k_{d,p}}, \\
\left\Vert p_{F}\right\Vert _{\infty } &\leq &K_{d,p}\left\Vert 1\right\Vert
_{W_{F}^{1,p}}^{1+k_{d,p}}
\end{eqnarray*}%
and%
\begin{equation*}
p_{F}(x)=\sum_{i=1}^{d}\E(\partial _{i}Q_{d}(F-x)\partial
_{i}^{F}1)=\sum_{i=1}^{d}\E(\partial _{i}Q_{d}(F-x)H_{i}(F;1)).
\end{equation*}

\textbf{B.}
For any positive $a<\frac{1}{d}-\frac{1}{p}$ one has
\begin{equation*}
p_{F }(x)\leq \Theta _{F}(\overline{p})\big(d+\left\Vert 1\right\Vert
_{W_{F}^{1,p}}\big)\,\P(F\in B_{2}(x))^{a}
\end{equation*}%
where $\overline{p}=1/(a+\frac{1}{p}).$
\end{theorem}

Now we give the representation formula and the estimates for the conditional expectation.

\begin{theorem}\label{th4-bis}
Suppose $1\in W_{F}^{1,p}$. Let $m\geq 1$ and $G\in W_{F}^{m,p}.$

\medskip

\textbf{A.}
We have $\mu_{F,G}(dx)=p_{F,G}(x)dx $ and
$$
\phi _{G}(x)=\E(G\mid F=x)=1_{\{p_{F}(x)>0\}}\frac{p_{F,G}(x)}{p_{F}(x)}
$$
with
$$
p_{F,G}(x)=\sum_{i=1}^{d}\E(\partial
_{i}Q_{d}(F-x)\partial _{i}^{F}G)=\sum_{i=1}^{d}\E(\partial
_{i}Q_{d}(F-x)H_{i}(F;G)).
$$

\textbf{B. }
We have $p_{F,G}\in W^{m,p}$ and $\phi _{G}\in \cap _{\delta >0}W^{m,p}(\{p_{F}>\delta \}).$
Moreover,
\begin{eqnarray*}
i)& &\left\Vert p_{F,G}\right\Vert _{\infty } \leq K_{d,p}\left\Vert
1\right\Vert _{W_{F}^{1,p}}^{k_{d,p}}\left\Vert G\right\Vert _{W_{F}^{1,p}}
\\
ii)& &\left\Vert p_{F,G}\right\Vert _{W^{m,p}} \leq
(2dK_{d,p})^{1-1/p}\left\Vert 1\right\Vert
_{W_{F}^{1,p}}^{k_{d,p}(1-1/p)}\left\Vert G\right\Vert _{W_{F}^{m,p}}.
\end{eqnarray*}%
We also have the representation formula
\begin{equation*}
\partial _{\alpha }p_{F,G}(x)=\sum_{i=1}^{d}\E(\partial
_{i}Q_{d}(F-x)\partial _{(\alpha ,i)}^{F}G)=\sum_{i=1}^{d}\E(\partial
_{i}Q_{d}(F-x)H_{(\alpha ,i)}(F;G))
\end{equation*}%
for any $\alpha $ with $0\leq
\left\vert \alpha \right\vert \leq m-1$. Furthermore,
$p_{F,G}\in C^{m-1}(\R^{d})$ and for any multi index $\alpha$ with $|\alpha|=k\leq m-2$,
$\partial_\alpha p_{F,G}$ is Lipschitz continuous with Lipschitz constant
$$
L_\alpha=d^2\|1\|_{W^{1,p}_F}^{k_d,p}\,\|G\|_{W^{k+2,p}_F}.
$$
And for any multi index $\alpha$ with $|\alpha|=m-1$,
$\partial_\alpha p_{F,G}$ is H\"older continuous of
exponent $1-d/p$ and H\"older constant
$$
L_\alpha=C_{d,p}\,\|p_{F,G}\|_{W^{m,p}},
$$
$C_{d,p}$ being dependent on $d$ and $p$ only.

\end{theorem}

Finally we give a stability property.

\begin{proposition}\label{stability}
Let $F_{n},G_{n},n\in\N$ be two sequences of random variables such that $%
(F_{n},G_{n})\rightarrow (F,G)$ in probability. Suppose that $G_{n}\in
W_{F_{n}}^{m,p}$ and $\sup_{n}(\left\Vert G_{n}\right\Vert
_{W_{F_{n}}^{m,p}}+\left\Vert F_{n}\right\Vert _{L^{p}})$ $<\infty $ for some $%
m\in \N.$ Then $G\in W_{F}^{m,p}$ and $\left\Vert G\right\Vert
_{W_{F}^{m,p}}\leq \sup_{n}\left\Vert G_{n}\right\Vert _{W_{F_{n}}^{m,p}}.$
\end{proposition}

\begin{remark}
Suppose that we are in the framework of Malliavin calculus and think that $F$
is a functional on the Wiener space which is non degenerated and
sufficiently smooth in Malliavin sense. And $G$ is another functional which
is sufficiently smooth in Malliavin sense. Then the Malliavin calculus
produces integration by parts formulas and so permits to prove that $G\in
W_{F}^{m,p}.$ But we may proceed in a different way: we start by taking a
sequence $F_{n},n\in \N$ of simple functionals such that $F_{n}\rightarrow F$
and a sequence $G_{n},n\in \N$ such that $G_{n}\rightarrow G$ and then we
may use standard finite dimensional integration by parts formulas in order to
prove that $G_{n}\in W_{F_{n}}^{m,p}.$ If we are able to check that $%
\sup_{n}\left\Vert G_{n}\right\Vert _{W_{F_{n}}^{m,p}}<\infty $ then using
the above stability property we conclude that $G\in W_{F}^{m,p}.$
\end{remark}

\textbf{Proof}. We denote $Q_{n}=(F_{n},G_{n},\partial _{\alpha
}^{F_{n}}G_{n},\left\vert \alpha \right\vert \leq m).$ Since $p\geq 2$ and $%
\sup_{n}(\left\Vert G_{n}\right\Vert _{W_{F_{n}}^{m,p}}$ $+\left\Vert
F_{n}\right\Vert _{L^{p}})<\infty $ it follows that the sequence $Q_{n},n\in
\N$ is bounded in $L^{2}$ and consequently weakly relative compact. Let $Q$
be a limit point. Using Mazur's theorem we construct a sequence of convex
combinations $\overline{Q}_{n}=\sum_{i=1}^{k_{n}}\lambda _{i}^{n}Q_{n+i},($%
with $\sum_{i=1}^{k_{n}}\lambda _{i}^{n}=1$ and $\lambda _{i}^{n}\geq 0)$
such that $\overline{Q}_{n}\rightarrow Q$ strongly in $L^{2}.$ And passing
to a subsequence we may assume that the convergence holds almost surely as
well. Since $(F_{n},G_{n})\rightarrow (F,G)$ in probability it follows that $%
Q=(F,G,\theta _{\alpha },\left\vert \alpha \right\vert \leq m).$ And using
the integration by parts formulas $IP_{\alpha ,p}(F_{n},G_{n})$ and the
almost sure convergence it is easy to see that $IP_{\alpha ,p}(F,G)$ holds
with $H_{\alpha }(F;G)=\theta _{\alpha }$ so $\theta _{\alpha }=\partial
_{\alpha }^{F}G.$\ Moreover using again the almost sure convergence and the
convex combinations one checks that $\left\Vert G\right\Vert
_{W_{F}^{m,p}}\leq \sup_{n}\left\Vert G\right\Vert _{W_{F_{n}}^{m,p}}.$ In
all the above arguments we have to use the Lebesgue dominated convergence
theorem so the almost sure convergence is not sufficient. But a
straightforward truncation argument which we do not develop here permits to
handle this difficulty. $\square $

\section{Functionals on the Wiener space}\label{sect-acc}

In this section we consider a probability space $(\Omega ,\cl F,\P)$ with a
Brownian motion $W=(W^{1},\ldots ,W^{n})$ and we use the Malliavin calculus in
order to obtain integration by parts formulas. We refer to D. Nualart \cite{bib:[N]}
for notation and basic results. We denote by $\D^{k,p}$ the space of the
random variables which are $k$ times differentiable in Malliavin sense in $%
L^{p}$ and for a multi-index $\alpha =(\alpha _{1},\ldots ,\alpha _{m})\in
\{1,\ldots ,n\}^{m}$ we denote by $D^{\alpha }F$ the Malliavin derivative of $F$
corresponding to the multi-index $\alpha .$
Moreover,
for any multi-index $\alpha$ with length  $\vert \alpha \vert=m$
we set
$$
\vert D^{(m)}F\vert ^{2}=\sum_{\vert \alpha \vert
=m}\vert D^{\alpha }F\vert ^{2}.
$$
We also consider the norms%
\begin{equation*}
\left\Vert F\right\Vert _{m,p}^{p}=\left\Vert F\right\Vert
_{p}^{p}+\sum_{k=1}^{m}\sum_{\left\vert \alpha \right\vert
=k}\E\Big(\Big(\int_{[0,\infty )^{k}}\left\vert D_{s_{1},\ldots ,s_{k}}^{\alpha
}F\right\vert ^{2}ds_{1}\ldots ds_{k}\Big)^{p/2}\Big).
\end{equation*}%
So $\D^{m,p}$ is the closure of the space of the simple functionals with
respect to the norm $\|\cdot\|_{m,p}.$ Moreover, for $F=(F^{1},\ldots ,F^{d}),F^{i}\in \D^{1,2},$
one denotes by $\sigma _{F}$ the Malliavin covariance matrix associated to $%
F:$%
\begin{equation*}
\sigma _{F}^{i,j}=\left\langle DF^{i},DF^{j}\right\rangle
=\sum_{k=1}^{n}\int_{0}^{\infty }D_{s}^{k}F^{i}D_{s}^{k}F^{j}ds,\quad
i,j=1,\ldots ,d.
\end{equation*}%
We will assume the non-degeneracy condition%
\begin{equation}
(\det \sigma _{F})^{-1}\in \cap _{p\in \N}L^{p}.  \label{diff1}
\end{equation}%
Under this assumption the matrix $\sigma _{F}$ is invertible and we denote
by $\widehat{\sigma }_{F}$ the inverse matrix. We also denote by $\delta $\
the divergence operator (Skorohod integral) and by $L$ the Ornstein
Uhlembeck operator and we recall that if $F\in \cap _{p\in \N}\D^{2,p}$ then $%
F\in Dom(L).$ The following proposition gives the classical integration by
parts formula from Malliavin calculus.

\begin{proposition}\label{mall-ibp}
i) Let $F=(F^{1},\ldots ,F^{d})$ with $F^{1},\ldots ,F^{d}\in \cap _{p\in \N}\D^{2,p}$
and $G\in \cap _{p\in \N}\D^{1,p}.$ Assume that (\ref{diff1}) holds. Then for
every function $f\in C_{b}^{1}(\R^{d})\rightarrow \R$ and every $i=1,\ldots ,d$
one has
\begin{equation}\label{diff2}
\begin{array}{l}
\E(\partial _{i}f(F)G) =-\E(f(F)H_{i}(F,G))\quad \mbox{with}   \\
\displaystyle
H_{i}(F,G) =-\sum_{j=1}^{d}\delta (G\widehat{\sigma }_{F}^{ji}DF^{j})=%
-\sum_{j=1}^{d}\Big(G\widehat{\sigma }_{F}^{ji}L(F^{j})+\left\langle DF^{j},D(%
\widehat{\sigma }_{F}^{ji}\times G)\right\rangle \Big)
\end{array}%
\end{equation}
and $H_{i}(F;G)\in \cap _{p\in \N}L^{p}.$

\smallskip

ii) Suppose that $F^{1},\ldots ,F^{d}\in \cap _{p\in \N}\D^{k+1,p}$ and $G\in \cap
_{p\in \N}\D^{k,p}$ for some $k\in \N.$ Then for every multi-index $\alpha
=(\alpha _{1},\ldots ,\alpha _{k})\in \{1,\ldots ,d\}^{k}$ one has
\begin{equation}
\E(\partial _{\alpha }f(F)G)=\E(f(F)H_{\alpha }(F,G))\quad with\quad H_{\alpha
}(F,G)=H_{\alpha _{k}}(F,H_{(\alpha _{1},\ldots ,\alpha _{k-1})}(F,G))
\label{diff4}
\end{equation}%
and $H_{\alpha }(F;G)\in \cap _{p\in \N}L^{p}.$
\end{proposition}

Notice that with the notation from the previous section we have
\begin{equation*}
\partial _{i}^{F}G=-\E(\delta (G(\widehat{\sigma }_{F}DF)^{i})\mid F).
\end{equation*}%
So Proposition \ref{mall-ibp} says that if (\ref{diff1}) holds and $%
F^{1},\ldots ,F^{d}\in \cap _{p\in \N}\D^{k+1,p}$ and $G\in \cap _{p\in \N}\D^{k,p}$
then $G\in \cap _{p\in \N}W_{F}^{k,p}.$

As an immediate consequence of Proposition \ref{mall-ibp} and Theorem \ref{th4}
we obtain the following result.

\begin{proposition}\label{mall-abc}
i) Let $F=(F^{1},\ldots ,F^{d})$ with $F^{1},\ldots ,F^{d}\in \cap _{p\in \N}\D^{2,p}.$
Assume that (\ref{diff1}) holds. Then the law of $F$ is absolutely
continuous with respect to the Lebesgue measure on $\R^{d}$ and the density $%
p_{F}$ can be represented as%
\begin{equation}
p_{F}(x)=\E\Big(\sum_{i=1}^{d}\partial _{i}Q_{d}(F-x)H_{i}(F;1)\Big).
\label{diff5}
\end{equation}%
Moreover $p_{F}$ is H\"{o}lder continuous of any exponent $q<1.$ And there
exists some universal constants $C_{d}$ and $p_{i},i=1,\ldots ,5$, depending on
$d$ such that
\begin{equation}
p_{F}(x)\leq C_{d}\E((\det \sigma _{F})^{-p_{1}})^{p_{2}}\left\Vert
F\right\Vert _{2,p_{4}}^{p_{3}}(P(\left\vert F-x\right\vert \leq
2))^{1/p_{5}}.  \label{diff5bis}
\end{equation}%
In particular, $\lim_{x\rightarrow \infty }\left\vert x\right\vert
^{p}p_{F}(x)=0$ for every $p\in \N.$

\smallskip

ii) Suppose that $F^{1},\ldots ,F^{d}\in \cap _{p\in \N}\D^{k+1,p}.$ Then $%
p_{F}\in C^{k-1}(\R^{d})$ and for every multi-index $\alpha $ with $%
\left\vert \alpha \right\vert \leq k$ one has%
\begin{equation}
\partial _{\alpha }p_{F}(x)=\E\Big(\sum_{i=1}^{d}\partial
_{i}Q_{d}(F-x)H_{(\alpha ,i)}(F;1)\Big).  \label{diff6}
\end{equation}%
Moreover, for $\left\vert \alpha \right\vert \leq k-1,$ $\partial _{\alpha
}p_{F}$\ is H\"{o}lder continuous of any exponent $q<1.$ And\ there exists some
universal constants $C_{d}$ and $p_{i},i=1,\ldots ,5$, depending on $d$ such
that
\begin{equation}
\left\vert \partial _{\alpha }p_{F}(x)\right\vert \leq C_{d}\E((\det \sigma
_{F})^{-p_{1}})^{p_{2}}\left\Vert F\right\Vert
_{k+1,p_{4}}^{p_{3}}(P(\left\vert F-x\right\vert \leq 2))^{1/p_{5}}.
\label{diff6bis}
\end{equation}%
In particular, if $F\in \cap _{p\in \N}L^{p}$ then $\lim_{x\rightarrow \infty
}\left\vert x\right\vert ^{p}\left\vert \partial _{\alpha
}p_{F}(x)\right\vert =0$ for every $p\in \N.$
\end{proposition}

\begin{remark}
The gain with respect to the classical result concerns the regularity (in
Malliavin sense) required for $F:$ recall that in the standard statement of
this criterion one needs that $F^{1},\ldots ,F^{d}\in \cap _{p\in \N}\D^{d+1,p}$
in order to obtain the existence of a continuous density and $%
F^{1},\ldots ,F^{d}\in \cap _{p\in \N}\D^{d+k+1,p}$ in order that the density is $%
k $ times differentiable. Moreover, notice that the estimate given in (\ref%
{diff5bis}) depends on the Sobolev norms of order two - and the same
estimates involve Sobolev norms of order $d+1$ if one applies the standard
criterion.
\end{remark}

\begin{remark}
The absolute continuity criterion of Bouleau and Hirsh \cite{bib:[B.H]} asserts
that if $F_{j}\in \D^{1,2},j=1,\ldots ,d$ and $\sigma _{F}\neq 0$ $a.s.$ then the
law of $F$ is absolutely continuous with respect to the Lebesgue measure.
The results presented in this section do not permit to prove this criterion
because we need at least one integration by parts and then we need that $%
F_{j}\in \D^{2,2},j=1,\ldots ,d.$ But if this stronger regularity assumption
holds, and moreover, if the stronger non degeneracy assumption (\ref{diff1})
holds as well, we obtain a density which is H\"{o}lder continuous and not
only measurable.
\end{remark}

\begin{remark}
The representation formula (\ref{diff5}) has been used by  Kohatsu Higa
and  Yasuda in \cite{bib:[K-H1]} and \cite{bib:[K-H2]} in order to provide numerical approximation
schemes for the density of the law of a diffusion process. Notice that $%
\E(\left\vert \partial _{i}Q_{d}(F-x)\right\vert ^{2})=\infty .$ Then
a direct use of the representation based on the Riesz transform leads to
approximation schemes of infinite variance which consequently are not
implementable by Monte Carlo methods. This is why they used a truncation
argument and gave an estimate of the error due to truncation. For this
estimate they used an old version of the present paper (namely \cite{bib:[B.C]}).
\end{remark}

Finally we give a result concerning the strict positivity set $U_F=\{p_F>0\}$.
We define the matrix field
$$
C^{ij}_F(x)=\E(H_i(F;1)|F=x)\E(H_j(F;1)|F=x)
$$
and the distance
$$
d_F(x,y)=\inf\Big\{\int_0^1\langle C_F(\varphi_t)\dot\varphi_t,\dot\varphi_t\rangle^{1/2}\,dt
\,;\, \varphi\in A^F_{x,y}\Big\}
$$
where $A^F_{x,y}=\{\varphi\in C^1([0,1],U_F)\,:\, \varphi_0=x,\varphi_1=y\}$.
Then,
\begin{proposition}\label{prop-conj}
Let $F=(F^{1},\ldots ,F^{d})$ with $F^{1},\ldots ,F^{d}\in \cap
_{p\in \N}\D^{3,p}$, and assume that (\ref{diff1}) holds. Then,
for any sequence $\{x_n\}_n\subset U_F$,
$\lim_{n\to\infty}p_F(x_n)=0$ if and only if
$\lim_{n\to\infty}d_F(x_n,x_1)=\infty$.

\end{proposition}

\textbf{Proof. } By recalling that Hirsch and Song
\cite{bib:[H.S]} proved that $U_F$ is a connected set, the
statement immediately follows by applying  parts $i)$ and $ii)$ of
Proposition \ref{mall-conj-prop}. $\square$

\smallskip

We give now the representation theorem for the conditional expectation.

\begin{proposition}\label{mall-cond-exp}
Let $F=(F^{1},\ldots ,F^{d})$ be such that $F^{1},\ldots ,F^{d}\in \cap _{p\in \N}\D^{2,p}$
and let  $G\in $ $\cap _{p\in \N}\D^{1,p}.$ Assume that (\ref{diff1}) holds. Then%
\begin{equation}\label{diff7}
\begin{array}{rcl}
\phi _{G}(x) &:=&\E(G\mid F=x)=1_{\{p_{F}>0\}}\displaystyle\frac{p_{F,G}(x)}{p_{F}(x)}%
\quad \mbox{with}  \\
p_{F,G}(x) &=&\E\Big(\sum_{i=1}^{d}\partial _{i}Q_{d}(F-x)H_{i}(F;G)\Big).
\end{array}
\end{equation}%
Moreover $\phi _{G}\in \cap _{\delta >0}\cap _{p\in \N}W^{1,p}(p_{F}>\delta )$
and it is locally H\"{o}lder continuous of any exponent $q<1$ on $\{p_{F}>0\}.$
And if $F^{1},\ldots ,F^{d}\in \cap _{p\in \N}\D^{k+2,p}$ and $G\in \cap _{p\in
\N}\D^{k+1,p}$ then $\phi _{G}\in C^{k}(\{p_{F}>0\}).$
\end{proposition}

\begin{remark}
The representation formula (\ref{diff7}) has been already obtained by
Malliavin and Thalmaier in \cite{bib:[M.T]} and was the starting point in our work. But
they need more regularity, namely $F^{1},\ldots ,F^{d}\in \cap _{p\in
\N}\D^{d+2,p.}.$ This is because they need to know that a bounded density
exists and they use the standard criterion for this.
\end{remark}

\end{document}